\documentclass[a4paper,twoside]{article}
\usepackage{a4}
\usepackage{amssymb}
\usepackage{amsmath}
\usepackage{upref}
\usepackage[active]{srcltx}
\usepackage[pagebackref,colorlinks,citecolor=blue,linkcolor=blue,urlcolor=blue]{hyperref}
\usepackage[dvipsnames]{color}
\allowdisplaybreaks[2] % To make it possible for displaybreaks
\newcount\minutes \newcount\hours
\hours=\time
\divide\hours 60
\minutes=\hours
\multiply\minutes -60
\advance\minutes \time
\newcommand{\klockan}{\the\hours:{\ifnum\minutes<10 0\fi}\the\minutes}
\newcommand{\tid}{\today\ \klockan}
\newcommand{\prtid}{\smash{\raise 10mm \hbox{\LaTeX ed \tid}}}
\renewcommand{\prtid}{}
\makeatletter
\def\sectionmark#1{} %\markboth{{\sectnr #1}}{{\sectnr #1}}} %Journal
\def\subsectionmark#1{}
\newcommand{\sectnr}{\ifnum \c@secnumdepth >\z@
                 \thesection.\hskip 1em\relax \fi}
\def\@evenhead{\footnotesize\rm\thepage\hfil\leftmark\hfil\llap{\prtid}}
\def\@oddhead{\footnotesize\rm\rlap{\prtid}\hfil\rightmark\hfil\thepage}
\def\tableofcontents{\section*{Contents} %\@mkboth{Contents}{Contents}} %Journal
 \@starttoc{toc}}
\makeatother
\makeatletter
\def\@biblabel#1{#1.}
\makeatother
\makeatletter
\let\Thebibliography=\thebibliography
\renewcommand{\thebibliography}[1]{\def\@mkboth##1##2{}\Thebibliography{#1}
\addcontentsline{toc}{section}{References}
\frenchspacing % Maybe not needed
\setlength{\@topsep}{0pt}% Delete if extra space before list
\setlength{\itemsep}{0pt}%
\setlength{\parskip}{0pt plus 2pt}%
}
\makeatother
\makeatletter
\def\mdots@{\mathinner.\nonscript\!.%
 \ifx\next,.\else\ifx\next;.\else\ifx\next..\else
 \nonscript\!\mathinner.\fi\fi\fi}
\let\ldots\mdots@
\let\cdots\mdots@
\let\dotso\mdots@
\let\dotsb\mdots@
\let\dotsm\mdots@
\let\dotsc\mdots@
\def\vdots{\vbox{\baselineskip2.8\p@ \lineskiplimit\z@
    \kern6\p@\hbox{.}\hbox{.}\hbox{.}\kern3\p@}}
\def\ddots{\mathinner{\mkern1mu\raise8.6\p@\vbox{\kern7\p@\hbox{.}}%
    \raise5.8\p@\hbox{.}\raise3\p@\hbox{.}\mkern1mu}}
\makeatother
\makeatletter
\let\Enumerate=\enumerate
\renewcommand{\enumerate}{\Enumerate%
\setlength{\@topsep}{0pt}% Delete if extra space before list
\setlength{\itemsep}{0pt}%
\setlength{\parskip}{0pt plus 1pt}%
\renewcommand{\theenumi}{\textup{(\alph{enumi})}}%
\renewcommand{\labelenumi}{\theenumi}%
}
\let\endEnumerate=\endenumerate
\renewcommand{\endenumerate}{\endEnumerate\unskip}
\makeatother
\newcommand{\addjustenumeratemargin}[1]{%
\setbox0\hbox{(a)} % Label to adjust width to 
\setbox1\hbox{#1} % Our labels here 
\addtolength{\leftmargini}{-\wd0}%
\addtolength{\leftmargini}{\wd1}%
}
\makeatletter
\def\@seccntformat#1{\csname the#1\endcsname.\quad}
\makeatother
\newcommand{\authortitle}[2]{\author{#1}\title{#2}\markboth{#1}{#2}}
\newcommand{\auth}[2]{{#1, #2.}}
\newcommand{\art}[6]{{\sc #1, \rm #2, \it #3\/ \bf #4 \rm (#5), \mbox{#6}.}}

\newcommand{\artprep}[3]{{\sc #1, \rm #2, \rm #3.}}
\newcommand{\artin}[3]{{\sc #1, \rm #2,  in #3.}}

\newcommand{\book}[3]{{\sc #1, \it #2, \rm #3.}}
\newcommand{\AND}{{\rm and }}
\RequirePackage{amsthm}
\newtheoremstyle{descriptive}%
  {\topsep}   %{\medskipamount}          % Space above
  {\topsep}   %  {\medskipamount}          % Space below
  {\rmfamily} % Body font
  {}          % Indent
  {\bfseries} % Head font
  {.}         % Punctuation after thm head
  { }         % Space after thm head
  {}          % Thm head spec(?)
\newtheoremstyle{propositional}%
  {\topsep}   %  {\medskipamount}          % Space above
  {\topsep}   %  {\medskipamount}          % Space below
  {\itshape}  % Body font
  {}          % Indent
  {\bfseries} % Head font
  {.}         % Punctuation after thm head
  { }         % Space after thm head
  {}          % Thm head spec(?)
\newtheoremstyle{remarkstyle}%
  {\topsep}   %  {\medskipamount}          % Space above
  {\topsep}   %  {\medskipamount}          % Space below
  {\rmfamily}  % Body font
  {}          % Indent
  {\itshape} % Head font
  {.}         % Punctuation after thm head
  { }         % Space after thm head
  {}          % Thm head spec(?)
\theoremstyle{propositional}
\newtheorem{thm}{Theorem}[section]
\newtheorem{prop}[thm]{Proposition}
\newtheorem{lem}[thm]{Lemma}
\newtheorem{cor}[thm]{Corollary}
\theoremstyle{descriptive}
\newtheorem{deff}[thm]{Definition}
\newtheorem{example}[thm]{Example}
\newtheorem{remark}[thm]{Remark}
\makeatletter
\renewenvironment{proof}[1][\proofname]{\par
  \pushQED{\qed}%
  \normalfont
  \trivlist
  \item[\hskip\labelsep
        \itshape
    #1\@addpunct{.}]\ignorespaces
}{%
  \popQED\endtrivlist\@endpefalse
}
\makeatother
\newcommand{\setm}{\setminus}
\renewcommand{\subsetneq}{\varsubsetneq}
\renewcommand{\supsetneq}{\varsupsetneq}
\renewcommand{\emptyset}{\varnothing}
\def\vint{\mathop{\mathchoice%
          {\setbox0\hbox{$\displaystyle\intop$}\kern 0.22\wd0%
           \vcenter{\hrule width 0.6\wd0}\kern -0.82\wd0}%
          {\setbox0\hbox{$\textstyle\intop$}\kern 0.2\wd0%
           \vcenter{\hrule width 0.6\wd0}\kern -0.8\wd0}%
          {\setbox0\hbox{$\scriptstyle\intop$}\kern 0.2\wd0%
           \vcenter{\hrule width 0.6\wd0}\kern -0.8\wd0}%
          {\setbox0\hbox{$\scriptscriptstyle\intop$}\kern 0.2\wd0%
           \vcenter{\hrule width 0.6\wd0}\kern -0.8\wd0}}%
          \mathopen{}\int}
\newcommand{\Cp}{{C_p}}
\DeclareMathOperator{\capp}{cap}
\newcommand{\cpDp}{\capp_{\Dp}}
\DeclareMathOperator{\dist}{dist}
\DeclareMathOperator{\spt}{supp}
\newcommand{\supp}{\spt}
\newcommand{\bdry}{\partial}
\newcommand{\bdy}{\bdry}
\newcommand{\loc}{_{\rm loc}}
\DeclareMathOperator{\Mod}{Mod}
\newcommand{\Modp}{\Mod_p}
\DeclareMathOperator{\para}{par}
{\catcode`p =12 \catcode`t =12 \gdef\eeaa#1pt{#1}}      % Get slantfactor
\def\accentadjtext#1{\setbox0\hbox{$#1$}\kern   % Convert it with height
                \expandafter\eeaa\the\fontdimen1\textfont1 \ht0 }
\def\accentadjscript#1{\setbox0\hbox{$#1$}\kern % Convert it with height
                \expandafter\eeaa\the\fontdimen1\scriptfont1 \ht0 }
\def\accentadjscriptscript#1{\setbox0\hbox{$#1$}\kern   % Convert it with height
                \expandafter\eeaa\the\fontdimen1\scriptscriptfont1 \ht0 }
\def\accentadjtextback#1{\setbox0\hbox{$#1$}\kern       % Convert it with height
                -\expandafter\eeaa\the\fontdimen1\textfont1 \ht0 }
\def\accentadjscriptback#1{\setbox0\hbox{$#1$}\kern     % Convert it with height
                -\expandafter\eeaa\the\fontdimen1\scriptfont1 \ht0 }
\def\accentadjscriptscriptback#1{\setbox0\hbox{$#1$}\kern % Convert it with height
                -\expandafter\eeaa\the\fontdimen1\scriptscriptfont1 \ht0 }
\def\itoverline#1{{\mathsurround0pt\mathchoice
        {\rlap{$\accentadjtext{\displaystyle #1}
                \accentadjtext{\vrule height1.593pt}
                \overline{\phantom{\displaystyle #1}
                \accentadjtextback{\displaystyle #1}}$}{#1}}
        {\rlap{$\accentadjtext{\textstyle #1}
                \accentadjtext{\vrule height1.593pt}
                \overline{\phantom{\textstyle #1}
                \accentadjtextback{\textstyle #1}}$}{#1}}
        {\rlap{$\accentadjscript{\scriptstyle #1}
                \accentadjscript{\vrule height1.593pt}
                \overline{\phantom{\scriptstyle #1}
                \accentadjscriptback{\scriptstyle #1}}$}{#1}}
        {\rlap{$\accentadjscriptscript{\scriptscriptstyle #1}
                \accentadjscriptscript{\vrule height1.593pt}
                \overline{\phantom{\scriptscriptstyle #1}
                \accentadjscriptscriptback{\scriptscriptstyle #1}}$}{#1}}}}
\newcommand{\al}{\alpha}
\newcommand{\alp}{\alpha}
\newcommand{\ga}{\gamma}
\newcommand{\la}{\lambda}
\newcommand{\om}{\omega}
\newcommand{\Om}{\Omega}
\newcommand{\clOm}{{\overline{\Om}}}
\newcommand{\clB}{\itoverline{B}}
\renewcommand{\phi}{\varphi}
\newcommand{\p}{{$p\mspace{1mu}$}}
\newcommand{\R}{\mathbf{R}}
\newcommand{\Sphere}{\mathbf{S}}
\newcommand{\wt}{\widetilde{w}}
\newcommand{\limplus}{{\mathchoice{\vcenter{\hbox{$\scriptstyle +$}}}
  {\vcenter{\hbox{$\scriptstyle +$}}}
  {\vcenter{\hbox{$\scriptscriptstyle +$}}}
  {\vcenter{\hbox{$\scriptscriptstyle +$}}}
}}
\newcommand{\Np}{N^{1,p}}
\newcommand{\Nploc}{N^{1,p}\loc}
\newcommand{\Dploc}{D^p\loc}
\newcommand{\Dp}{D^p}
\newcommand{\Ga}{\Gamma}
\newcommand{\Galoc}{\Gamma\loc}
\newcommand{\Lploc}{L^p\loc}
\newcommand{\La}{\Lambda}
\newcommand{\Xhat}{{\widehat{X}}}
\newcommand{\lahat}{{\hat{\la}}}
\newcommand{\uhat}{{\hat{u}}}
\newcommand{\vhat}{{\hat{v}}}
\newcommand{\ghat}{{\hat{g}}}
\newcommand{\CpX}{{C_p^X}}
\newcommand{\Nhat}{\widehat{N}}
\newcommand{\muhat}{\hat{\mu}}
\newcommand{\phihat}{\widehat{\phi}}
\numberwithin{equation}{section}
\newenvironment{ack}{\medskip{\it Acknowledgement.}}{}
\newcommand{\veq}{{\scriptstyle\|}}
\newcommand{\vsubset}{\cup}
\usepackage{mathtools}
\renewcommand{\veq}{\mathrel{\vcenter{\hbox{\rotatebox{90}{$=$}}}}}
\renewcommand{\vsubset}{\mathrel{\vcenter{\hbox{\rotatebox{90}{$\subset$}}}}}
\begin{document}

\authortitle{Anders Bj\"orn, Jana Bj\"orn
    and Nageswari Shanmugalingam}
{Classification of metric measure spaces and their ends using
    \p-harmonic functions} 
\author{
Anders Bj\"orn \\
\it\small Department of Mathematics, Link\"oping University, SE-581 83 Link\"oping, Sweden\\
\it \small anders.bjorn@liu.se, ORCID\/\textup{:} 0000-0002-9677-8321
\\
\\
Jana Bj\"orn \\
\it\small Department of Mathematics, Link\"oping University, SE-581 83 Link\"oping, Sweden\\
\it \small jana.bjorn@liu.se, ORCID\/\textup{:} 0000-0002-1238-6751
\\
\\
Nageswari Shanmugalingam
\\
\it \small  Department of Mathematical Sciences, University of Cincinnati,
P.O.\ Box 210025,\\
\it \small   Cincinnati, OH 45221-0025, U.S.A.\/{\rm ;}
\it \small shanmun@uc.edu, ORCID\/\textup{:} 0000-0002-2891-5064
}

\date{Preliminary version, \today}
\date{}
\maketitle

\begin{center}
Dedicated to the memory of Seppo Rickman.
\end{center}

\noindent{\small
{\bf Abstract}. 
  By seeing whether
a Liouville type theorem holds for positive, bounded, and/or finite energy
\p-harmonic and \p-quasiharmonic functions,
we classify proper 
metric spaces equipped
with a locally doubling measure supporting a local \p-Poincar\'e inequality.
Similar classifications have earlier been obtained for Riemann surfaces
  and Riemannian manifolds.

We study the inclusions between these classes of metric measure spaces, and
their relationship
  to the \p-hyperbolicity 
  of the metric space and its ends.
In particular, we
characterize spaces that carry nonconstant \p-harmonic functions with finite energy
as spaces having at least two well-separated \p-hyperbolic sequences.
  We also show that 
  every such space $X$ has a function
  $f \notin L^p(X) + \R $ with finite \p-energy.
}

\medskip

\noindent
    {\small \emph{Key words and phrases}:
classification of metric measure spaces,
doubling measure,
end at infinity,
finite energy,
\p-hyperbolic sequence,
Liouville theorem,
\p-harmonic function, 
Poincar\'e inequality,
\p-parabolic,
quasiharmonic function,
quasiminimizer.}

\medskip


\section{Introduction}

The classical Liouville theorem states that every bounded holomorphic
 function in the whole plane is constant.
A similar statement is true for harmonic and \p-harmonic functions in
$\R^n$, $1<p<\infty$.

In the 1960s, Riemann surfaces were classified according
to existence of global analytic or harmonic functions in various classes
(bounded, positive and finite-energy), which culminated in the
1970 monograph by Sario and~Nakai~\cite{SaNa70}.
Together with Wang and Chung, they extended this classification
to Riemannian manifolds in the monograph~\cite{SNWC} from 1977.
Holopainen~\cite{HoThesis} extended this classification
further to \p-harmonic functions on Riemannian manifolds in 1990,
see also Kilpel\"ainen~\cite[Theorem~1.8]{Kilp89} for some
similar results for Euclidean domains.
Subsequently, in the
1990s, first-order analysis
on metric spaces began to be studied and it has since seen
a growing interest.
Our main aim is to obtain a similar 
classification of metric spaces as in the monographs mentioned above.
We refer to later sections for the definitions.

\medskip

\noindent
\emph{Throughout the paper, except for
Sections~\ref{sect-prelim} and~\ref{sect-noncomp},
we assume that $1<p<\infty$ and
that $X$ is an unbounded proper connected metric space
equipped with a locally doubling measure $\mu$ supporting
a local \p-Poincar\'e inequality.}

\medskip

\begin{deff}\label{def:Op}
  We say that  $X$ belongs to the \emph{Liouville type class}
\addjustenumeratemargin{$O^p_{HBD}$}% To adjust leftmargin in enumerate
\begin{enumerate}
\item[$O^p_{HP}$] if every \emph{positive} \p-harmonic function on $X$ is constant\/\textup{;}
\item[$O^p_{HB}$] if every \emph{bounded} \p-harmonic function on $X$ is constant\/\textup{;}
\item[$O^p_{HD}$] if every \p-harmonic function on $X$ with \emph{finite energy} is constant\/\textup{;}
\item[$O^p_{HBD}$] if every \emph{bounded} \p-harmonic function on $X$
  with \emph{finite energy} is constant. 
\end{enumerate}
The corresponding classes
$O^p_{QP}$, $O^p_{QB}$, $O^p_{QD}$ and $O^p_{QBD}$
for   quasiharmonic functions (where the dependence on $p$ is implicit)
are defined similarly.
Moreover, we say that  $X \in O^p_{\para}$ if $X$ is \emph{\p-parabolic}
in the sense of Definition~\ref{def-hyp-end}.
\end{deff}

Our classification result can be summarized as follows.

\begin{thm} \label{thm-intro-class}
We have the following inclusions\/\textup{:}
\begin{equation*}
  \begin{matrix}
    O^p_{HP} & \subsetneq & O^p_{HB} & \subset & O^p_{HBD} & =  & O^p_{HD} &
    \supsetneq &  O^p_{\para} \\
   \vsubset & & \vsubset && \veq && \veq \\
  O^p_{QP} & \subsetneq &  O^p_{QB} & \subset & O^p_{QBD} & =  & O^p_{QD}. & 
      \end{matrix}
\end{equation*}
Moreover, $O^p_{QB} \setm O^p_{HP}$,
$O^p_{QP} \setm O^p_{\para}$
  and $ O^2_{HBD}\setm O^2_{HB}$ are nonempty.
\end{thm}

Some of these inclusions are of course trivial.
In the setting of orientable Riemannian manifolds,
it was shown in
Sario--Nakai--Wang--Chung~\cite{SNWC} (for $p=2$) and 
Holopainen~\cite{HoThesis} (for general $p$)
  that
\begin{equation} \label{eq-O2}
O^p_{\para}\subsetneq O^p_{HP}\subsetneq O^p_{HB}\subset O^p_{HBD}=O^p_{HD}
  \quad \text{and} \quad 
  O^2_{HB}\subsetneq O^2_{HBD}.
\end{equation}
(Whenever we discuss manifolds we implicitly assume
  that they are connected and have dimension $\ge 2$.)

A class of functions called
``quasiharmonic'' was also considered 
in~\cite{SNWC}.
However, those functions  
  are 
solutions to $\Delta u = 1$, while our quasiharmonic functions are continuous
quasiminimizers of the \p-energy.
Such quasiminimizers
 were introduced in Giaquinta--Giusti~\cite{GG1}, \cite{GG2}
  as a unified treatment of variational inequalities, elliptic
partial differential equations and quasiregular mappings, 
see \cite{BBS-Liouville} and~\cite{JBqmin15}
for further discussion and references.

Since Riemann, planar Euclidean domains have been classified using conformal mappings:
two planar domains belong to the same category if there is a conformal mapping
between them.
One of the motivations for 
studying the classes in~\eqref{eq-O2} is that some of them
are conformally invariant
on Riemann surfaces, when $p=2$.
Consequently, two conformally equivalent Riemann
surfaces either both belong to such a 
class or neither belongs to
that class.

For higher-dimensional Euclidean domains and $p\ne2$,
conformal mappings are too rigid, and instead
quasiconformal or quasisymmetric mappings are used.
For $n$-dimensional Riemannian
manifolds, $n$-harmonicity and $n$-parabolicity are conformal
invariants.

The theory of quasiconformal mappings was extended to
metric measure spaces in
Heinonen--Koskela~\cite{HeKo95}, \cite{HeKo98},
see also
Heinonen--Koskela--Shan\-mu\-ga\-lin\-gam--Ty\-son~\cite[Section~9]{HKST01}.
Quasiconformal mappings do not preserve harmonic or \p-harmonic functions,
but they do preserve quasiharmonic functions (with $p=Q$)
in proper connected
spaces with a
 uniformly locally Ahlfors $Q$-regular measure supporting
 a uniformly local $Q$-Poincar\'e inequality,
see Korte--Marola--Shan\-mu\-ga\-lin\-gam~\cite[Theorem~4.1]{KMS12}
and also Heinonen--Kilpel\"ainen--Martio~\cite[Corollary~4.7]{HeKi92}.
Quasiconformal mappings between such spaces therefore
  preserve the classes $O^Q_{QB}$ and $O^Q_{QP}$,
  and, by
\cite[Theorem~9.10]{HKST01},
  also $O^Q_{QD}$ and $O^Q_{QBD}$.
Hence it is
natural to include quasiharmonic Liouville type classes in our study.
The existence of nonconstant global quasiharmonic functions on one (but
  not the other) space therefore
gives a convenient way of checking whether two metric measure spaces can be
quasiconformally equivalent.
As far as we know,
  even for $p=2$ and in the setting of Riemann surfaces,
   it is not known whether the classes $O^2_{HP}$ and $O^2_{HB}$
  are quasiconformally invariant,
see~Sario--Nakai~\cite[p.~7]{SaNa70}.
On the other hand, it was noted already therein
that $O^2_{HD}$ is quasiconformally
  invariant in that setting.

For complete 
Riemannian manifolds,
the case $p=2$ is also related 
to the
Brownian motion:
$2$-parabolicity is equivalent
to the fact that almost surely the Brownian motion
starting from a compact set $K$ will intersect each neighborhood of
$K$ infinitely often,
  see Grigor$'$yan~\cite[Theorem~5.1]{Grig99}.
Thus the classification of metric measure spaces
as in Theorem~\ref{thm-intro-class} has roots 
in the theory of Brownian motion,
in complex dynamics (see~\cite[Theorem~0.1]{NR}),
and in the study of quasiconformal maps.

A natural way of distinguishing between different spaces
and manifolds is through their ends at infinity.
For instance, (unweighted) $\R^n$ has one end if $n \ge 2$
and this end is \p-hyperbolic if and only if $1<p < n$.
For $n=1$, $\R$ has two ends which are both \p-parabolic.
An end, or a space, is \p-hyperbolic if it is not \p-parabolic,
see Definition~\ref{def-hyp-end}.

We show that if $X$ has two \p-hyperbolic ends, then $X \notin O^p_{HBD}$.
The converse is not true
as explained in Example~\ref{ex-Tree-hyp-but-parEnds}, but using
the new concept of \emph{\p-hyperbolic sequences} we are able to
give the following characterization.

\begin{thm} \label{thm-intro-hyp-seq}
$X \notin O^p_{HBD}$ if and only if 
there are two disjoint \p-hyperbolic sequences
$\{F_n\}_{n=1}^\infty$ and $\{G_n\}_{n=1}^\infty$
which are well-separated in the sense that 
the \p-modulus of the family $\Ga(F_1,G_1)$ of all curves from $F_1$ to $G_1$
satisfies
\[
\Modp(\Ga(F_1,G_1)) <\infty.
\]
In this case, $X$ is also \p-hyperbolic,
i.e.\ $O^p_{\para} \subset O^p_{HBD}$.

In particular, $X \notin O^p_{HBD}$ if $X$ has two \p-hyperbolic ends.
\end{thm}

As (unweighted) $\R^n \in O^p_{HP} \subset O^p_{HBD}$ for all $1<p<\infty$ and $n \ge 1$,
but
is \p-parabolic only for $p \ge n$, 
we see  that
$O^p_{\para} \subsetneq O^p_{HBD}$ (cf.\ Theorem~\ref{thm-intro-class}).

For $p=2$, similar characterizations of the bounded and finite-energy Liouville
theorems
(i.e.\ of $X \in O^2_{HB}$ resp.\ $X \in O^2_{HD}$)
by means of well-separated massive
and/or hyperbolic sets were obtained for Riemannian manifolds, see
Grigor$'$yan~\cite[Proposition~1 and Theorem~2]{Grig87}, 
\cite[Theorem~13.10\,(b)]{Grig99} and the references therein.
  In the setting of Gromov hyperbolic spaces with  uniformly local
assumptions (of doubling and \p-Poincar\'e inequality), the validity
of the finite-energy Liouville theorem
for \p-harmonic functions (i.e.\ $X \in O^p_{HD}$)
was characterized using
uniformization in
Bj\"orn--Bj\"orn--Shan\-mu\-ga\-lin\-gam~\cite[Theorem~10.5]{BBSunifPI}.
See Remark~\ref{rmk-unifPI} for how our results in this paper
  improves upon that.

Hyperbolic sequences
can be seen as subsets of the 
hyperbolic parts of the ``boundary at $\infty$'' of the metric space $X$. 
For simply
connected complete Riemannian manifolds $M$ of negative sectional
curvature, such a ``boundary at $\infty$'', $M(\infty)$, was introduced by
Eberlein--O'Neill~\cite{EON} and identified with the ``sphere at
$\infty$''.
If $M$, in addition, has negatively pinched sectional curvature 
$-b^2\le K\le -a^2 <0$,
then it is possible to solve the asymptotic Dirichlet problem with
any continuous boundary data on the sphere at infinity.
This follows from Choi~\cite[Theorems~4.5 and~4.7]{C} and
Anderson~\cite{And83} for $p=2$,
and has been generalized to $p>1$  by 
Pansu~\cite{Pan} and Holopainen~\cite[Theorem~2.1]{Ho02}, see
also the discussion in \cite[p.~3394]{Ho02}.
These existence results imply that $M\notin O^p_{HB}$, but do
not address the existence of \p-harmonic functions with finite energy.
On Gromov hyperbolic spaces, the above
solvability of the Dirichlet problem at infinity was 
deduced
in Holopainen--Lang--V\"ah\"akangas~\cite[Theorem~6.2]{HLV} for $p>1$,
under various additional assumptions.

Choi~\cite[Definition~5.1]{C} considers ends on a finitely connected complete
$2$-dimensional Riemannian manifold
with sectional curvature $K\le-a^2 <0$ and
shows that if the surface is orientable with
$\int K=-\infty$, then
it carries many nonconstant bounded harmonic functions, see
\cite[Theorem~5.13 and Corollary~5.14]{C}. 
The Dirichlet problem for \p-harmonic functions in unbounded domains
with ends towards infinity was solved in
Bj\"orn--Bj\"orn--Li~\cite[Theorems~6.6, 7.6 and~7.7]{BBLi} in the
setting of Ahlfors $Q$-regular spaces under certain assumptions on
$Q$, $p$ and the measure. 
The notions of parabolic and hyperbolic ends
have also been used in the study of some other
partial differential equations,
see for instance Korolkov--Losev~\cite{KL} for the case of
the stationary Schr\"odinger equation.

Under global assumptions, the Liouville theorem (Theorem~\ref{thm-Liouville})
for positive quasiharmonic functions on metric spaces was
obtained in Kinnunen--Shan\-mu\-ga\-lin\-gam~\cite{KiSh01}.   
In~\cite{BBS-Liouville}, we proved
the so-called finite-energy Liouville theorem 
for noncomplete spaces
with global assumptions under various additional assumptions.
We can now deduce
the Liouville theorem for finite-energy 
quasiharmonic functions
without those additional assumptions,
as a direct consequence of our identity $O^p_{QBD}=O^p_{QD}$
(and Theorem~\ref{thm-Liouville})
provided that $X$ is complete, see
Corollary~\ref{cor-finite-energy-Liouville}.
Moreover, using 
tools from Bj\"orn--Bj\"orn~\cite{BBnoncomp}
we are able to lift this also to noncomplete spaces,
see Theorem~\ref{thm-Dp-Liouville-noncomplete}.
The weighted real line and Example~\ref{ex-weighted-R2} show that
  the finite-energy Liouville theorem
  fails if the global assumptions are relaxed to uniformly local ones.

The following theorem shows that
the Liouville type class $O^p_{HD}$ is
  related to the question of whether every function with finite
  energy on $X$ can be written as a global Sobolev function plus a
  constant, i.e.\ whether $D^p(X)=\Np(X)+\R$.
this is the case for $p=2$, then 
the classical theory of Dirichlet forms and
the associated spectral decomposition
can be extended to the Dirichlet space $D^{2}(X)$ of functions with
finite energy, where the associated Dirichlet form is in terms of the 
Cheeger differential structure as in
Koskela--Rajala--Shan\-mu\-ga\-lin\-gam~\cite{KRS}.

\begin{thm} \label{thm:main2}
If $X\notin O^p_{HD}$, then $\Dp(X)\ne N^{1,p}(X)+\R$.
\end{thm}

Example~\ref{ex:OpHD-but-not-hom} shows that the converse fails.
However, if $X\in O^p_{HD}$ and $X$ supports a global
$(p,p)$-Sobolev inequality (in addition to our standing assumptions),
then $\Dp(X)=N^{1,p}(X)+\R$,
see Proposition~\ref{prop-OpHD-andHom}.

The rest of this paper is structured as follows. Definitions of the
concepts related to the function spaces studied in this paper are
given in Section~\ref{sect-prelim},
and the concepts regarding \p-harmonicity and related useful tools are given
in Section~\ref{sect-qmin}. 
Section~\ref{sect-ends} is devoted to the definitions of \p-hyperbolic ends
and \p-hyperbolic sequences in metric measure spaces
and a brief discussion of them.
In Section~\ref{sect-existence} we prove the existence of nonconstant
\p-harmonic functions with finite energy
under the assumption that the metric measure space has at least
two distinct \p-hyperbolic sequences.
We follow this up by a discussion of classification of metric measure spaces
in Section~\ref{sec:Classification}.
In this section we also provide the proofs of
Theorems~\ref{thm-intro-class} and~\ref{thm-intro-hyp-seq}.
The third main theorem of the paper,
Theorem~\ref{thm:main2}, is proved in Section~\ref{sect-Dirichlet}.
In that section the converse of 
Theorem~\ref{thm:main2} is also discussed.
Section~\ref{sec-ex} is devoted to providing examples that 
  illustrate the sharpness of the results given
  in the paper.
    Example~\ref{ex-R-hyp-par} is also essential when
    deducing most of the noninclusions in Theorem~\ref{thm-intro-class}.
Finally, 
Section~\ref{sect-noncomp} provides an 
  extension of the finite-energy Liouville theorem to the noncomplete
  setting.

\begin{ack}
Parts of this research project were conducted during
2017 and 2018 when N.~S. was a guest professor at Link\"oping University, partially funded
by the  Knut  and  Alice  Wallenberg  Foundation,  and  during  the  parts  of 2019 when 
A. B. and J. B. were Taft Scholars at the University of Cincinnati. 
The authors would like to thank these institutions for their kind
support and hospitality. A.~B. and J.~B. were partially supported by
the Swedish Research Council grants
2016-03424 and 2020-04011 resp.\ 621-2014-3974 and 2018-04106.
N.~S. was partially supported by the National Science Foundation (U.S.A.)
grants DMS-1500440 and DMS-1800161.
\end{ack}

\section{Preliminaries}
\label{sect-prelim} 

We assume throughout the paper that
$X$ is a metric space equipped
with a metric $d$ and a positive complete  Borel  measure $\mu$ 
such that $0<\mu(B)<\infty$ for all balls $B \subset X$.
In this section we also assume that $1 \le p< \infty$.
For proofs of the facts stated in this section 
we refer the reader to Bj\"orn--Bj\"orn~\cite{BBbook} and
Heinonen--Koskela--Shan\-mu\-ga\-lin\-gam--Tyson~\cite{HKST}.

A notion critical to this paper is that of \p-modulus of families of curves in $X$. 
A \emph{curve} is a continuous mapping from an interval.
We will only consider locally rectifiable curves,
  and they can always be parameterized by their arc length $ds$.

\begin{deff} \label{deff-modulus}
  Let $\Gamma$ be a family of locally rectifiable  curves in $X$. 
The \emph{\p-modulus} of $\Gamma$ is the number
\[
\Modp(\Gamma):=\inf_\rho \int_X\rho^p\, d\mu,
\]
where the infimum is taken over all nonnegative 
Borel functions $\rho$ on $X$ such that 
$\int_\ga\rho\, ds\ge 1$ for each $\ga\in\Ga$.
\end{deff}

From now on, unless otherwise said, all our curves will be nonconstant, compact
and \emph{rectifiable}, i.e.\ of finite length.
We follow Heinonen and Koskela~\cite{HeKo98} in introducing
upper gradients (in~\cite{HeKo98} they are referred to 
as very weak gradients).

\begin{deff} \label{deff-ug}
A nonnegative Borel function $g$ on $X$ is an \emph{upper gradient} 
of an extended real-valued function $u$ on $X$ if for all curves  
$\gamma: [0,l_{\gamma}] \to X$,
\begin{equation} \label{ug-cond}
        |u(\gamma(0)) - u(\gamma(l_{\gamma}))| \le \int_{\gamma} g\,ds,
\end{equation}
where we follow the convention that the left-hand side is $\infty$ 
whenever at least one of the terms therein is infinite.
If $g$ is a nonnegative measurable function on $X$
and if (\ref{ug-cond}) holds for \p-almost every curve, 
then $g$ is a \emph{\p-weak upper gradient} of~$u$. 
A property holds for \emph{\p-almost every curve}
if it fails only for a  curve family 
with zero \p-modulus.
\end{deff}

The notion of \p-weak upper gradients was introduced in
Koskela--MacManus~\cite{KoMc}. It was also shown therein
that if $g \in \Lploc(X)$ is a \p-weak upper gradient of $u$,
then one can find a sequence $\{g_j\}_{j=1}^\infty$
of upper gradients of $f$ such that $\|g_j-g\|_{L^p(X)} \to 0$.

If $u$ has an upper gradient in $\Lploc(X)$, then
it has a \emph{minimal \p-weak upper gradient} $g_u \in \Lploc(X)$
in the sense that $g_u \le g$ a.e.
for every \p-weak upper gradient $g \in \Lploc(X)$ of $u$,
see Shan\-mu\-ga\-lin\-gam~\cite{Sh-harm}.
The minimal \p-weak upper gradient is well defined
up to a set of measure zero in the cone of nonnegative functions in $\Lploc(X)$.
Moreover, 
$g_u=g_v$ a.e.\ in the set $\{x \in X : u(x)=v(x)\}$,
in particular $g_{\min\{u,c\}}=g_u \chi_{\{u < c\}}$ for $c \in \R$.
Note also that a modification of an upper gradient on a Borel set of
measure zero need not
yield an upper gradient, but a modification of a \p-weak upper gradient on a 
set of measure zero still yields a \p-weak upper gradient.

Following Shan\-mu\-ga\-lin\-gam~\cite{Sh-rev}, 
we define a version of Sobolev spaces on $X$. 

\begin{deff} \label{deff-Np}
For a measurable function $u:X\to [-\infty,\infty]$, let 
\[
        \|u\|_{\Np(X)} = \biggl( \int_X |u|^p \, d\mu 
                + \inf_g  \int_X g^p \, d\mu \biggr)^{1/p},
\]
where the infimum is taken over all upper gradients $g$ of $u$.
The \emph{pre-Newtonian space} on $X$ is 
\[
        \Np (X) = \{u: \|u\|_{\Np(X)} <\infty \}.
\]
\end{deff}

The \emph{Newtonian space} $\Np(X)/{\sim}$,
where  $f \sim h$ if and only if $\|f-h\|_{\Np(X)}=0$,
is a Banach space and a lattice, 
see~\cite{Sh-rev}. 
We are also interested in the homogeneous version of
  Sobolev spaces. 
The \emph{Dirichlet space} $\Dp(X)$ is the collection of all
measurable functions on $X$ 
that have an upper gradient in $L^p(X)$.

We say  that $u \in\Nploc(X)$ if
for every $x \in X$ there exists $r_x$ such that 
$u \in\Np(B(x,r_x))$.
The local spaces $L^p\loc(X)$ and $\Dp\loc(X)$ are defined similarly.
Note that if $X$ supports a local \p-Poincar\'e inequality (as in
Definition~\ref{def-local} below) then it follows by truncations and
Fatou's lemma that $\Nploc(X)=\Dp\loc(X)$.

In this paper we assume that functions 
in the above function spaces $\Nploc(X)$, $\Dp\loc(X)$
are defined everywhere (with values in $[-\infty,\infty]$),
not just up to an equivalence class in the corresponding function space.

For a measurable set $E\subset X$, the 
space $\Np(E)$ is defined by
considering $(E,d|_E,\mu|_E)$ as a metric space in its own right.
The spaces $\Nploc(E)$, $L^p(E)$, $\Lploc(E)$,
$\Dp(E)$ and $\Dploc(E)$ are defined similarly.

\begin{deff} \label{deff-Cp}
The (Sobolev) \emph{capacity} of a set $E\subset X$  is the number 
\begin{equation*} 
   \CpX(E)=\Cp(E) =\inf_u    \|u\|_{\Np(X)}^p,
\end{equation*}
where the infimum is taken over all $u\in \Np (X) $ such that $u=1$ on $E$.
\end{deff}

A property is said to hold \emph{quasieverywhere} (q.e.) 
if the set of all points in $X$ at which the property
fails has $\Cp$-capacity zero. 
The capacity is the correct gauge 
for distinguishing between two Newtonian functions. 
If $u \in \Np(X)$, then $u \sim v$ if and only if $u=v$ q.e.
Moreover, if $u,v \in \Nploc(X)$ and $u= v$ a.e., then $u=v$ q.e.

\begin{deff}
The (Dirichlet) capacity of the pair $(E,F)$ of disjoint sets in $X$ is
\[
\cpDp(E,F)=\int_X g_u^p\, d\mu,
\] 
where the infimum is taken over all functions $u\in \Dp(X)$ with $u\ge 1$ on $E$
and $u\le 0$ on $F$. 
\end{deff}

The following equality was proved for compact sets in
  Kallunki--Shan\-mu\-ga\-lin\-gam~\cite{KaS01}. Since we need it for general
  closed sets, we provide a short proof.
  Here and later we let $\Ga(E,F)$ 
be the collection of all 
curves in $X$ with one end point in $E$ and the other in $F$.
  
\begin{lem}   \label{lem-Mod=capDp}
  Let $E$ and $F$ be disjoint closed subsets of $X$.
Then 
\[
\Modp(\Ga(E,F)) = \cpDp(E,F).
\]
\end{lem}

\begin{proof}
Let $v\in\Dp(X)$ be admissible for $\cpDp(E,F)$.
Then every upper gradient $g$ of $v$ is admissible for $\Modp(\Ga(E,F))$ 
and hence 
\[
\Modp(\Ga(E,F)) \le \int_X g^p\,d\mu.
\]
Taking infimum over all upper gradients $g$ of $v$ and then taking
infimum over all $v$ admissible for $\cpDp(E,F)$ proves one inequality
in the lemma.

Conversely, let $\rho\in L^p(X)$ be admissible for
$\Modp(\Ga(E,F))$ and consider the function
\[
u(x):= \min\biggl\{1, \inf_\ga \int_\ga\rho\,ds\biggr\},
\]
with the infimum taken over all 
curves (including constant curves) $\ga$ connecting
$x$ to~$F$. By Bj\"orn--Bj\"orn--Shan\-mu\-ga\-lin\-gam~\cite[Lemma~3.1]{BBS5},
$u$ has $\rho$ as an upper gradient,
$u=0$ on $F$ and $u=1$ on $E$.
Since $\rho\in L^p(X)$, we infer from
Rogovin--Rogovin--J\"arvenp\"a\"a--J\"arvenp\"a\"a--Shan\-mu\-ga\-lin\-gam~\cite[Corollary~1.10]{JJRRS}
that $u$ is measurable and thus $u\in \Dp(X)$. It follows that 
\[
\cpDp(E,F) \le \int_X \rho^p\,d\mu,
\]
and taking infimum over all $\rho\in L^p(X)$ admissible for 
$\Modp(\Ga(E,F))$ concludes the proof.
\end{proof}

As in Bj\"orn--Bj\"orn~\cite{BBsemilocal}, we
define the following local 
versions of the notions of doubling measures and Poincar\'e inequality.

\begin{deff} \label{def-local}
We say that the measure $\mu$ is \emph{doubling within a ball $B_0$}
if there is a \emph{doubling constant} $C>0$ (depending on $B_0$) 
such that for all balls $B=B(x,r):=\{y\in X: d(y,x)<r\} \subset B_0$,
\[
\mu(2B)\le C \mu(B),
\]
where $\lambda B=B(x,\lambda r)$.  

Similarly,  the \emph{\p-Poincar\'e inequality holds within a ball $B_0$}
if there are constants $C>0$ and $\lambda \ge 1$
(both depending on $B_0$)
such that for all balls $B\subset B_0$,
all integrable functions $u$ on $\la B$, and all upper gradients $g$ of $u$
in $\la B$, 
\begin{equation}  \label{eq-PI-on-B}
        \vint_{B} |u-u_B| \,d\mu
        \le C r_B \biggl( \vint_{\lambda B} g^{p} \,d\mu \biggr)^{1/p},
\end{equation}
where $u_B:=\vint_B u \,d\mu := \int_B u\, d\mu/\mu(B)$
and $r_B$ is the radius of $B$. 

Each of these properties is
called \emph{local} if 
for every $x \in X$ there is some $r>0$
(depending on $x$) such that the property 
holds within $B(x,r)$. The property is
called \emph{uniformly local} if $r$, $C$ and $\la$ are
independent of $x$. If it holds within every ball $B(x_0,r_0)$ in $X$ with
$C$ and $\la$ 
independent of $x_0$ and $r_0$, then
it is 
called \emph{global}.
\end{deff}

\section{Quasiharmonic and \texorpdfstring{\p}{p}-harmonic functions}
\label{sect-qmin}

\emph{From now on, except for Section~\ref{sect-noncomp},
  we assume that $X$ is an unbounded proper connected metric space.
We also assume that $1<p<\infty$, that $\mu$ is locally doubling 
  and supports a local \p-Poincar\'e inequality, and that 
  $\Om \subset X$ is an open set.}
  
\medskip 

A metric space $X$ is \emph{proper} if every closed bounded set in $X$ is compact.
It follows that $X$ is complete. Moreover, Proposition~1.2 and Theorem~1.3
in Bj\"orn--Bj\"orn~\cite{BBsemilocal} imply that under the above assumptions,
the doubling property and \p-Poincar\'e inequality actually hold
within every ball in $X$.

\begin{deff} \label{def-qmin}
A function $u \in \Np\loc(\Om)$ is a \emph{quasiminimizer} in $\Om$ if
there exists $Q_u\ge1$ such that
\begin{equation} \label{eq-deff-qmin}
     \int_{\phi \ne 0} g^p_u \,d\mu \le Q_u  \int_{\phi \ne 0} g^p_{u+\phi} \,d\mu
\end{equation}
for all $\phi \in \Np_0(\Om)$, where
\[
\Np_0(\Om):=
  \{\phi|_{\Om} : \phi \in \Np(X) \text{ and }
        \phi=0 \text{ on } X \setm \Om\}.
\]
A \emph{quasiharmonic function} is a continuous quasiminimizer.

If $Q_u=1$ in~\eqref{eq-deff-qmin}, then $u$ is a \emph{minimizer},
and if it is in addition continuous, then it is a
\emph{\p-harmonic function}.
\end{deff}

Functions from $\Np_0(\Om)$ can be extended by zero in $X\setm \Om$ and we
will regard them in that sense if needed.

Note that 
the property of being a quasiminimizer 
depends on the index $p$
even though we have refrained from making that explicit in the notation.
The integrals in \eqref{eq-deff-qmin} can be infinite but then
they are infinite simultaneously.
Under our assumptions, locally Lipschitz functions are dense in
$\Np\loc(\Om)$, see \cite[Theorem~8.4]{BBsemilocal}.
It therefore follows from
Bj\"orn--Bj\"orn--Shan\-mu\-ga\-lin\-gam~\cite[Theorem~5.7]{BBS5}
  (or \cite[Theorem~5.45]{BBbook})
that Lipschitz functions 
with compact support in $\Om$ are dense in $\Np_0(\Om)$.
Hence, the definition of quasiminimizers can equivalently be based on such
compactly supported Lipschitz test functions.
The
  integration in~\eqref{eq-deff-qmin} can moreover equivalently be over 
  $\supp\phi$ rather than the set where $\phi\ne0$, see
 Bj\"orn~\cite[Proposition~3.2]{ABkellogg}.
Note also that $\Np_0(X)=\Np(X)$, which has consequences for globally 
defined  quasiminimizers on $X$.

Any quasiminimizer can be modified on a set of capacity zero
so that it becomes locally H\"older continuous.
This follows from the results in 
Kinnunen--Shan\-mu\-ga\-lin\-gam~\cite[p.~417]{KiSh01}.
The assumptions therein are different from ours, but
see Bj\"orn--Bj\"orn~\cite[Theorem~10.2 and the discussion
    around it]{BBsemilocal}
for how those results apply under the local assumptions considered here.
Such a continuous representative  is called a
\emph{quasiharmonic} function or, for $Q_u=1$, a \emph{\p-harmonic function}.

The Liouville theorem given below follows from the Harnack inequality proved 
in~\cite[Corollary~7.3]{KiSh01} 
or Bj\"orn--Marola~\cite[Corollary~9.4]{BMarola}.

\begin{thm} \label{thm-Liouville}
 Assume that $\mu$ is globally doubling and supports a global
  \p-Poincar\'e inequality.
  If $u$ is a positive quasiharmonic function on $X$, then it is constant.
  In particular, $X \in O^p_{QP}\subset O^p_{HP}$.
\end{thm}

The following lemma
will be convenient when proving Theorem~\ref{thm:main2}.

\begin{lem}  \label{lem-Liouville-Np}
  If $u\in\Np(X)$ is quasiharmonic on $X$, then it is constant.
\end{lem}

\begin{proof}
Let $u$ be quasiharmonic on $X$.
If $u\in\Np(X) = \Np_0(X)$, then testing~\eqref{eq-deff-qmin} with $-u
  \in \Np_0(X)$ yields 
\[
     \int_{u \ne 0} g^p_u \,d\mu \le Q_u  \int_{u \ne 0} g^p_{u-u} \,d\mu=0.
\]
This together with the local \p-Poincar\'e inequality shows that
$u$ is locally a.e.-constant,
and as $u$ is continuous and $X$ connected, $u$ is constant.
\end{proof}

The following lemma about convergence of \p-harmonic functions
is
a useful tool.

\begin{lem}   \label{lem-reflex-conv}
Let $\Om_j$ be open sets such that $\Om_j\subset\Om_{j+1}$, $j=1,2,\ldots$\,, and
$X=\bigcup_{j=1}^\infty\Om_j$.
Assume that $u_j\in\Dp(X)$ is \p-harmonic in $\Om_j$ and that there
is a constant $M$ such that for all $j=1,2,\ldots$\,,
\[
|u_j|\le M \text{ in } X \quad \text{and} \quad \|g_{u_j}\|_{L^p(X)} \le M.
\]
Then there are\/ \textup{(}finite\/\textup{)} convex combinations
\[
\uhat_j=\sum_{k=j}^{N_j}\tilde{\la}_{k,j}u_k, 
\quad \text{with} \ \tilde{\la}_{j,k}\ge0 \ \text{and} \
  \sum_{k=j}^{N_j} \tilde{\la}_{j,k}=1,
\]
of the sequence $\{u_j\}_{j=1}^\infty$, 
which converge
locally uniformly in $X$ to a function $u\in\Dp(X)$
that is \p-harmonic in $X$,
satisfies $|u|\le M$ and moreover 
\begin{equation}   \label{eq-grad-uhat-to-gu-Lp}
\|g_{\uhat_j}-g_u\|_{L^p(X)}\to 0, \quad \text{as } j\to\infty.  
\end{equation}
\end{lem}

\begin{proof}
Theorem~5.4 in~\cite{BBsemilocal} implies that for every ball
$B_0\subset X$, there is some $1\le q<p$ such that a $q$-Poincar\'e
inequality holds within this ball in the sense of Definition~\ref{def-local}.
This better Poincar\'e inequality allows us to apply the 
continuity and convergence results for \p-harmonic functions from
Kinnunen--Shan\-mu\-ga\-lin\-gam~\cite{KiSh01} and Shan\-mu\-ga\-lin\-gam~\cite{Sh03},
see also the discussion in~\cite[Section~10]{BBsemilocal}.

More precisely, by 
\cite[Proposition~3.3 and Theorem~5.2]{KiSh01}
and the fact that $|u_j|\le M$ on $X$, the (tail of the)
sequence $\{u_j\}_{j=1}^\infty$ is equi(H\"older)-continuous on every ball
in $X$, see also \cite[Theorem~8.14]{BBbook}.
Thus an appeal to the Ascoli theorem and the Harnack convergence principle 
(\cite[Theorem~1.2]{Sh03} or \cite[Theorem~9.37]{BBbook}),
together with a Cantor diagonalization argument, 
yields a subsequence, also denoted $\{u_j\}_{j=1}^\infty$,
that converges uniformly on balls in $X$ to a function $u$ that is \p-harmonic in $X$.
Note that $u\in\Np(B)$ for every ball $B$ and that $|u| \le M$ on $X$.
It remains to prove~\eqref{eq-grad-uhat-to-gu-Lp}.

Since the sequence $\{g_{u_j}\}_{j=1}^\infty$ is bounded in $L^p(X)$,
we can use the reflexivity of $L^p(X)$ to extract
a subsequence, still denoted $\{g_{u_j}\}_{j=1}^\infty$, that
converges weakly to a nonnegative function $g\in L^p(X)$. 
Mazur's lemma (applied iteratively to the subsequences $\{g_{u_j}\}_{j=k}^\infty$) 
then provides us with a sequence of convex combinations 
\[
g_k =\sum_{j=k}^{N(k)}\la_{j,k} g_{u_j}, 
\quad \text{with} \ \la_{j,k}\ge0 \ \text{and} \ \sum_{j=k}^{N(k)}\la_{j,k}=1,
\]
such that $\|g_k-g\|_{L^p(X)}\le 2^{-k}$.
Let $\ghat = g + \sum_{k=1}^\infty |g_k-g|$. 
Then $\ghat\in L^p(X)$ and $g_k\le \ghat$ 
in $X$ for all $k=1,2,\ldots$\,.

Note that the functions $g_k$ are \p-weak upper gradients of the
corresponding convex combinations
\[
v_k =\sum_{j=k}^{N(k)}\la_{j,k} u_j.
\]
Hence $g_{v_k}\le g_k$ a.e.\ in $X$ and
\begin{equation}  \label{eq-bound-g-vj}
  \|g_{v_k}\|_{L^p(X)} \le \|g_k\|_{L^p(X)}
  \le \sum_{j=k}^{N(k)}\la_{j,k} \|g_{u_k}\|_{L^p(X)} \le M.
\end{equation}

Next, choose an increasing sequence of balls $B_j$, so that $X=\bigcup_{j=1}^\infty B_j$.
The sequence $\{v_k\}_{k=1}^\infty$ satisfies $|v_k|\le M$ on $X$.
In view of~\eqref{eq-bound-g-vj}, it is therefore bounded in
$N^{1,p}(\clB_j)$ for every $j=1,2,\ldots$\,.
Since $\clB_j$ is a complete doubling metric space
by Bj\"orn--Bj\"orn~\cite[Propositions~1.2 and~3.4]{BBsemilocal},
it follows from 
Ambrosio--Colombo--Di Marino~\cite[Corollary~41]{AmbCD} that
the Newtonian space $\Np(\clB_j)/{\sim}$ is reflexive
(where $f \sim h$ if and only if $\|f-h\|_{\Np(X)}=0$).
Thus, using weakly converging subsequences and Mazur's lemma again
(this time for subsequences in $\Np(\clB_j)$), 
for each $j=1,2,\ldots$ we can find a further convex combination 
\[
\uhat_j =\sum_{k=j}^{\Nhat(j)}\lahat_{j,k} v_k, 
\quad \text{with} \ \lahat_{j,k}\ge0 \ \text{and} \ \sum_{k=j}^{\Nhat(j)}\lahat_{j,k}=1,
\]
such that $\|\uhat_j-u\|_{\Np(\clB_j)} \le 2^{-j}$.
In particular, $\|g_{\uhat_j-u}\|_{L^p(B_j)} \le 2^{-j}$. As
\[
g_u \le g_{\uhat_j}+g_{u-\uhat_j} \quad \text{and} \quad g_{\uhat_j} \le g_u+g_{\uhat_j-u},
\]
we consequently have
$\|g_{\uhat_j}-g_u\|_{L^p(B_j)} \le 2^{-j}$.
In particular, $g_{\uhat_j}\to g_u$ in $L^p(B)$
for each ball $B$ and a.e.\ in $X$, as $j\to\infty$.

Note that the sequence $\{\uhat_j\}_{j=1}^\infty$ (being a
  convex combination of locally uniformly converging functions)
  also converges locally uniformly in $X$ to $u$.

Now, since $g_{v_k} \le g_k \le \ghat$ a.e.\ in $X$, we conclude that also
\[
g_{\uhat_j} \le \sum_{k=j}^{\Nhat(j)}\lahat_{j,k} g_{v_k} \le \ghat \in L^p(X).
\]
The Lebesgue dominated convergence theorem therefore implies that
$g_{\uhat_j}\to g_u$ in $L^p(X)$, which concludes the proof.
\end{proof}

We will use solutions of the Dirichlet problem and more precisely
so-called \p-harmonic extensions, which we define
next, following Hansevi~\cite[Definition~4.6]{hansevi1}. 

\begin{deff}  \label{def-harm-ext}
Assume that $\Cp(X \setm \Om)>0$.
Let $f \in \Dp(X)$.
Then the \emph{\p-harmonic extension}
$H_\Om f$ of $f$ in $\Om$ is the unique \p-harmonic function in $\Om$
such that $f-H_\Om f \in \Dp_0(\Om)$,
where
\[
\Dp_0(\Om):=
  \{\phi|_{\Om} : \phi \in \Dp(X) \text{ and }
        \phi=0 \text{ on } X \setm \Om\}.
\]
We also let $H_\Om f=f$ on $X \setm \Om$ to get a globally
  defined function when needed.
\end{deff}

The \p-harmonic extension exists  
and is unique by Hansevi~\cite[Theorem~4.4]{hansevi1}.
If $\Om$ is bounded and $f \in \Np(X)$, then the definition
of $H_\Om f$ coincides with other definitions in the literature, such
as in Shan\-mu\-ga\-lin\-gam~\cite[Theorem~5.6]{Sh-harm},
Bj\"orn--Bj\"orn--Shan\-mu\-ga\-lin\-gam~\cite[Definition~3.3]{BBS}
and \cite[Definition~8.31]{BBbook}.
The existence, uniqueness and other properties of $H_\Om f$ in bounded
sets were obtained in these references.

The following relation between $\cpDp$ and harmonic extensions
  in unbounded sets will be useful.

\begin{prop} \label{prop-cpDp-Hf}
Let $F_0$ and $F_1$ be two disjoint closed sets with
$\cpDp(F_0,F_1)<\infty$.
Then there is $f \in \Dp(X)$ such that $f=j$ on $F_j$, $j=0,1$.
Moreover for any such $f$,
\begin{equation} \label{eq-cpDp-Hf}
\cpDp(F_0,F_1)  = \int_X g_{H_\Om f}^p \, d\mu,
\quad \text{where }
\Om=X \setm (F_0 \cup F_1).
\end{equation}
\end{prop}  

\begin{proof}
As $\cpDp(F_0,F_1)<\infty$, the existence of such a function $f$ is immediate.
The definition of the harmonic extension 
in~\cite[Definition~4.6]{hansevi1}
shows that $H_\Om f$ solves 
the minimization problem in the definition
of $\cpDp(F_0,F_1)$, i.e.\ it satisfies \eqref{eq-cpDp-Hf}.
\end{proof}

\begin{remark}   \label{rem-local-ass-Hansevi}
Since we consider \p-harmonic functions on unbounded sets in this
paper, results from Hansevi~\cite{hansevi1}, \cite{hansevi2} will be of primary
importance here.
We therefore comment on how the assumptions therein compare with ours.

In~\cite{hansevi1}, $X$ is assumed to be proper, connected and supporting a global
$(p,p)$-Poincar\'e inequality (with an averaged
$L^p$-norm also on the left-hand side of \eqref{eq-PI-on-B}).
However, the only use of the Poincar\'e inequality
in the existence theorem~\cite[Theorem~3.4]{hansevi1},
the comparison principle~\cite[Lemma~3.6]{hansevi1} and the
convergence theorem for obstacle problems~\cite[Theorem~3.2]{hansevi2} is through
Maz$'$ya's inequality on a sequence of balls
(on p.~98 and again on p.~102) with no need of uniform control of the constants.
Therefore, it is enough to require that $\mu$ supports a 
$(p,p)$-Poincar\'e inequality on all sufficiently large balls,
see the proof of Maz$'$ya's inequality in \cite[Theorem~5.53]{BBbook}.
Under our standing assumptions, such a  
$(p,p)$-Poincar\'e inequality on balls (with constants depending
on the ball)
follows from Bj\"orn--Bj\"orn~\cite[Theorems~1.3 and~5.1]{BBsemilocal}.

The inner regularity results in~\cite[Theorem~4.4]{hansevi1} and the
tools from Kinnunen--Shan\-mu\-ga\-lin\-gam~\cite{KiSh01}, Shan\-mu\-ga\-lin\-gam~\cite{Sh03},
Kinnunen--Martio~\cite{KiMa02} and Bj\"orn--Bj\"orn~\cite{BBbook},
used in~\cite{hansevi1} and \cite{hansevi2}
are of local nature and
therefore hold under our local assumptions. 
Note that since $X$ is assumed to be connected and proper, the local
doubling property and Poincar\'e inequality self-improve so that they
actually hold within every ball $B_0\subset X$ (with constants
depending on $B_0$), which is enough for such local regularity results,
see the discussion in~\cite[Section~10]{BBsemilocal}.

In particular, the resolutivity and uniqueness results for Perron
solutions with continuous boundary data on unbounded \p-parabolic sets
from~\cite[Section~7]{hansevi2} are available under our assumptions
and will be used later.
\end{remark}

\section{Hyperbolic ends and hyperbolic sequences}
\label{sect-ends}

The theory of ends was originally developed to study
the classification of Riemann surfaces,
as in Sario--Nakai~\cite{SaNa70}.
Heuristically, for us an end represents a point of $X$ at $\infty$. 
For example, 
if $X$ is homeomorphic to
  $\Sphere^{1}\times\R$,
then in our sense it has two
ends.
However, note that
if the metric on $X$ is such that at least one of the ends is hyperbolic and rotationally invariant,
then from a geometric group theoretic point of view this end contains a copy of
$\Sphere^1$. In this paper 
we 
still consider this end as one point at $\infty$.

\begin{deff}\label{def:ends}
We say that a sequence $\{F_n\}_{n=1}^\infty$ is a \emph{chain} at $\infty$ of $X$ (called a chain of $X$ for simplicity)
if there is a point $x_0\in X$ and a strictly increasing sequence of radii
$R_n\to\infty$ such that $F_n$ is a component of
$X\setminus B(x_0,R_n)$ and $F_{n+1}\subset F_n$. 

Two chains $\{F_n\}_{n=1}^\infty$ and $\{G_n\}_{n=1}^\infty$
at $\infty$ are said to be \emph{equivalent} if for each positive
integer $k$ there are  
$n_k$ and $m_k$ such that $F_{n_k}\subset G_k$ and $G_{m_k}\subset F_k$. 
This equivalence relationship partitions the class of
all chains of $X$ into pairwise disjoint equivalence classes, 
called \emph{ends} of $X$.
\end{deff}

From Kline--Lindquist--Shan\-mu\-ga\-lin\-gam~\cite[Lemma~5.11]{KLS} 
we know that the choice of $x_0$ does not play a central role in the construction
of ends.
Traditionally, an end of a manifold or a metric space
$X$ is a sequence $\{F_n\}_{n=1}^\infty$ of
connected sets that are components of complements of compact subsets
$K_n \subset X$ such that $F_{n+1}\subset F_n$
for each $n$ and $X=\bigcup_{n=1}^\infty K_n$, see e.g.\ 
Choi~\cite[Definition~5.1]{C}.
For us it is more convenient
to have ends made up of closed sets. 
Given our assumption that $X$ is proper, replacing $F_n$
with its closure merely gives an equivalent notion of ends.

The papers  Grigor$'$yan~\cite{Grig94}, \cite{Grig99} and
Holopainen~\cite{Ho99}
used different definitions of ``ends'', sufficient for their
purposes.
Since in this paper we will be discussing the possibility of a metric space 
having more than one end
and even infinitely many ends,
we need the precise terminology here. 

The terminology we follow is adapted from 
Adamowicz--Bj\"orn--Bj\"orn--Shan\-mu\-ga\-lin\-gam~\cite{ABBS}
and Estep~\cite{EstThesis}.
The equivalence class that contains a chain $\{F_n\}_{n=1}^\infty$ was
denoted $[F_n]$ in~\cite{EstThesis} and~\cite{KLS}. 
However, it was shown in~\cite{EstThesis} that if $\{F_n\}_{n=1}^\infty$ and 
$\{G_n\}_{n=1}^\infty$ are two chains at $\infty$ such that 
for each $k$ there is a positive integer $m_k$ with $F_{m_k}\subset G_k$, then the 
two chains are equivalent. 
Therefore in discussing an end, it suffices to
discuss a chain at $\infty$ that represents the end. 
Hence from now on, we will also call a chain an end.

Recall the
definitions of $\Modp$ and $\Ga(E,F)$ from Section~\ref{sect-prelim}.
Given an end $\{F_n\}_{n=1}^\infty$ and $E\subset X$,  
note that
\[
\Gamma(E, F_{n+1}) \subset \Gamma(E, F_n).
\]
As in Holopainen~\cite{Ho99} and Holopainen--Koskela~\cite{HK01} 
we give the following definitions.
On metric measure spaces,
the study of \p-parabolicity and \p-hyperbolicity
began in~\cite{HK01}
and Holopainen--Shan\-mu\-ga\-lin\-gam~\cite{HS02}.

\begin{deff}   \label{def-hyp-end}
The end $\{F_n\}_{n=1}^\infty$ is \emph{\p-hyperbolic} if 
\[
\lim_{n\to\infty}\Modp(\Gamma(\itoverline{B(x_0,1)}, F_n))>0.
\]
The space $X$ is \emph{\p-hyperbolic} if
\[
\lim_{n\to\infty}\Modp(\Gamma(\itoverline{B(x_0,1)}, X\setm B(x_0,n)))>0.
\]
We say that an end or $X$ is \emph{\p-parabolic} if it is not \p-hyperbolic.  
\end{deff}

It follows directly from the definition that if $X$ has a
\p-hyperbolic end, then $X$ is a \p-hyperbolic space.
Conversely, if $X$ is a \p-hyperbolic space with 
finitely many ends, then it has a \p-hyperbolic end.
Example~\ref{ex-Tree-hyp-but-parEnds} below shows
that
there is a \p-hyperbolic space with infinitely many \p-parabolic
ends but no \p-hyperbolic end.

\begin{remark}\label{rem:end-paths}
Since the metric measure space $X$ is proper, we know 
from Shan\-mu\-ga\-lin\-gam~\cite[Theorem~4.2]{Sh-Fractals} that if $\{F_n\}_{n=1}^\infty$ is an end of~$X$, 
then it is \p-hyperbolic if and only if 
\[
\Modp (\Galoc(\clB,\{F_n\}_{n=1}^\infty)) >0,
\]
where $B=B(x_0,1)$ and
$\Galoc(\clB,\{F_n\}_{n=1}^\infty)$ is the
collection of all locally rectifiable curves $\ga$
starting in $\clB$ and intersecting $F_n$
for each  $n=1,2,\ldots$\,.
Note that
\[
\Galoc(\clB,\{F_n\}_{n=1}^\infty) =\bigcap_{n=1}^\infty\Galoc(\clB, F_n)
\quad \text{and} \quad
\Modp(\Galoc(\clB, F_n))=\Modp(\Gamma(\clB,F_n)),
\]
where $\Galoc(\clB, F_n)$
is the collection of all locally 
rectifiable curves in $X$ starting in $\clB$ and intersecting $F_n$.
\end{remark}

Euclidean spaces $\R^n$, $n \ge 2$,
have exactly one end, and this end is \p-parabolic if and only
  if $p\ge n$.
Parabolicity can in many situations be characterized by volume
  growth conditions, see 
\cite[Theorem~5.5]{BBLehIntGreen}, \cite[Proposition~3.4]{CHS01}, 
\cite[Theorems~7.3 and~14.6]{Grig99}, \cite[Section~4]{Ho99}
and~\cite[Theorem~1.7]{HK01}.

Our aim in this paper is to investigate
when a metric measure space
  carries nonconstant \p-(quasi)harmonic functions.
For functions with finite energy, this property turns out to be
closely  related to the following notion, which
extends the concepts given in Definition~\ref{def-hyp-end}.

\begin{deff}  \label{def-hyp-seq}
Let $x_0 \in X$.
  A sequence $\{F_n\}_{n=1}^\infty$ is a \emph{\p-hyperbolic sequence} if
  it is a decreasing sequence of nonempty closed sets such that\/\textup{:}
  \begin{enumerate}
  \item \label{l-i}
    for each $r>0$ there is $n>0$ such that $B(x_0,r) \cap F_n = \emptyset$\textup{;}
  \item
\begin{equation} \label{eq-def-hyp-seq}
\lim_{n\to\infty}\Modp(\Ga(\itoverline{B(x_0,1)}, F_n))>0.
\end{equation}
  \end{enumerate}
\end{deff}

Since $X$ is proper, it follows that
\ref{l-i} is equivalent to $\bigcap_{n=1}^\infty F_n=\emptyset$.
This equivalence need not hold in nonproper spaces.
It follows directly from the definitions that
every \p-hyperbolic end forms a \p-hyperbolic sequence,
and that the existence of a \p-hyperbolic sequence implies that
$X$ is
  \p-hyperbolic.
The following lemma shows that in Definitions~\ref{def-hyp-end} and~\ref{def-hyp-seq},
the ball $\itoverline{B(x_0,1)}$ can equivalently be replaced by any compact set
with positive capacity, see also Holopainen--Shan\-mu\-ga\-lin\-gam~\cite[Proof of Lemma~3.5]{HS02}.

\begin{lem}   \label{lem-cap-K1-K2}
Let $\{F_n\}_{n=1}^\infty$ be a decreasing sequence of closed sets in $X$
satisfying condition~\ref{l-i} of Definition~\ref{def-hyp-seq}. 
Also let $K_1$ and $K_2$ be compact sets with $\Cp(K_j)>0$, $j=1,2$.
Then
\[
\lim_{n\to\infty}\cpDp(K_1,F_n) = 0 \quad \text{if and only if} \quad
\lim_{n\to\infty}\cpDp(K_2,F_n) = 0.
\]
\end{lem}

\begin{proof}
By replacing $K_1$ resp.\ $K_2$ by $K_1 \cup K_2$ we see that
we may assume, without loss of generality, that $K_1 \subset K_2$.
It follows from \cite[Proposition~4.8 and Lemma~4.10]{BBsemilocal} that
there is a ball $B$ such
that $K_2$ is contained in a rectifiably pathconnected component $G$ of $B$.
Note that as $X$ is locally doubling and supports a local Poincar\'e inequality,
it is locally quasiconvex, and hence $G$ is an open set.

As $K_1$ is compact,
there is $f \in \Np(X)$ such that $f \equiv 0$ in $X\setm G$
and $f \equiv 1$ on $K_1$.
Let $v=H_{G \setm K_1} f$
be the \p-harmonic extension of $f$ in $G\setm K_1$
as in Definition~\ref{def-harm-ext}.
We start by showing that $m:=\inf_{K_2} v >0$.
If not, then the strong  maximum principle
(see~\cite{KiSh01} or~\cite[Theorem~8.13]{BBbook},
together with \cite[Section~10]{BBsemilocal}) shows that
$v\equiv 0$ in $G \setm K_1$.
Moreover, $v \equiv 0$ in $X \setm G$, and
since $v \equiv 1$
in $K_1$, we would have that $v=\chi_{K_1} \in \Np(X)$
and $g_v=0$ a.e.
  The \p-Poincar\'e inequality on $B$ then implies that $v$ is constant a.e.\
  (and thus q.e.) in $B$, which contradicts $\Cp(K_1)>0$.
  Note 
from~\cite[Theorem~1.3]{BBsemilocal} that our assumptions on $X$ imply
the validity of a \p-Poincar\'e inequality on arbitrary balls
(with constants depending on the ball).
Thus $m>0$.

Let $n$ be 
large enough so that $F_n\cap B =\emptyset$.
Let $u=H_{X \setm (K_1 \cup F_n)} f$ be the \p-harmonic extension of $f$
in $X \setm (K_1 \cup F_n)$.
By Proposition~\ref{prop-cpDp-Hf}, 
\[
\int_X g_u^p\,d\mu = \cpDp(K_1,F_n).
\]
The comparison principle (see for example Hansevi~\cite[Lemma~3.6]{hansevi1})
 implies that $u\ge v$ in $B$.
Hence $u\ge m$ in $K_2$,
and thus $u/m$ is admissible for $\cpDp(K_2,F_n)$.
Therefore
\[
m^p \cpDp(K_2,F_n) \le \int_X g_u^p\,d\mu 
         = \cpDp(K_1,F_n) \le \cpDp(K_2,F_n).
\]
Since $m$ is independent of $n$,         
letting $n\to\infty$ concludes the proof.
\end{proof}

  \section{Existence of nonconstant
    \texorpdfstring{\p}{p}-harmonic functions with finite energy} \label{sect-existence}

\begin{thm} \label{thm-2-hyp-seq}
  Assume that there are two disjoint
  \p-hyperbolic sequences   $\{F_n\}_{n=1}^\infty$ and $\{G_n\}_{n=1}^\infty$
  such that
$\Modp(\Ga(F_1,G_1))<\infty$.
  Then $X$ supports a 
  nonconstant bounded \p-harmonic function with finite energy,
    i.e.\ $X \notin O^p_{HBD}$.
\end{thm}

Observe that $\Modp(\Ga(F_1,G_1))=\cpDp(F_1,G_1)$, by Lemma~\ref{lem-Mod=capDp}.
Before proving Theorem~\ref{thm-2-hyp-seq}, we first show how
it implies the following important corollary.
At the same time, Example~\ref{ex-weighted-R2} shows that
Theorem~\ref{thm-2-hyp-seq} can be used also when $X$ only has one
end, and that end is \p-hyperbolic.
The converse of Theorem~\ref{thm-2-hyp-seq} is proved in
  Proposition~\ref{prop-p-harm-imp-2-hyp-seq} below.

\begin{cor}
\label{cor:mainly-here}
If $X$ has at least two \p-hyperbolic ends,  then $X$ supports a 
nonconstant bounded \p-harmonic function with finite energy,
    i.e.\ $X \notin O^p_{HBD}$.
  \end{cor}

\begin{proof}
A \p-hyperbolic end is automatically a \p-hyperbolic sequence.
Denote the two \p-hyperbolic ends as
$\{F_n\}_{n=1}^\infty$ and $\{G_n\}_{n=1}^\infty$, with $F_1 \cap G_1 = \emptyset$.
We may also assume that they are created as in Definition~\ref{def:ends}
with the same strictly increasing sequence of radii $\{R_n\}_{n=1}^\infty$.
Testing $\Modp(\Ga(F_2,G_2))$ with
\[
\rho=
\frac{\chi_{B(x_0,R_2)}}{R_2 - R_1}
\]
then shows that $\Modp(\Ga(F_2,G_2)) < \infty$.
Hence (after shifting indices), 
the corollary follows from Theorem~\ref{thm-2-hyp-seq}.
\end{proof}

To prove Theorem~\ref{thm-2-hyp-seq} we first need to understand connectivity 
properties of curves between two \p-hyperbolic sequences, or
equivalently the relationship between the corresponding capacities.
This will be the content of the following two lemmas.

\begin{lem}          \label{lem:hyp-seq-imp-Ga}
Suppose that $\{F_n\}_{n=1}^\infty$ and $\{G_n\}_{n=1}^\infty$ are two disjoint \p-hyperbolic 
sequences in $X$.
Then
\begin{equation*} 
 \lim_{n\to\infty}\cpDp(F_n,G_n)>0.
\end{equation*}
\end{lem}

Recall that $\cpDp(F_n,G_n)=\Modp(\Ga(F_n,G_n))$.

\begin{proof}
Because of the monotonicity of $\cpDp$
and choosing a subsequence if necessary,
without loss of generality
we may assume that $\cpDp(F_1,G_1)<\infty$. 
Let $B=B(x_0,1)$.
By Theorem~10.2 in Bj\"orn--Bj\"orn~\cite{BBsemilocal}
there are
positive constants $C$ and $\La$ such that 
the weak Harnack inequality
\begin{equation}   \label{eq-weak-Harnack-B_1}
\vint_{B} v\,d\mu \le C \inf_{B} v
\end{equation}
holds for all nonnegative \p-harmonic functions $v$ in
$\La B$.
From the definition of \p-hyperbolic sequences, we can find $N$ such that
$\La B\cap (F_N \cup G_N)=\emptyset$.
Let $n \ge N$ be fixed but arbitrary.
As $\cpDp(F_n,G_n)<\infty$, there is $f \in \Dp(X)$ with
$f=0$ on $F_n$ and  $f=1$ on $G_n$.
It follows from \eqref{eq-def-hyp-seq} and Lemma~\ref{lem-Mod=capDp}
that $\Cp(F_n)>0$.
By Proposition~\ref{prop-cpDp-Hf}
\[
\cpDp(F_n,G_n)=\int_X g_{u}^p\, d\mu,
\]
where $u=H_{X \setm (F_n \cup G_n)} f$ is the \p-harmonic extension
of $f$ in $X \setm (F_n \cup G_n)$.
Let
\[
m= \inf_{B} u \quad \text{and} \quad M= \sup_{B} u.
\]
We 
distinguish two cases. 
If $\vint_{B} u\,d\mu\ge\tfrac12$, then the weak Harnack 
inequality~\eqref{eq-weak-Harnack-B_1} implies that $2Cm\ge 1$
and hence any   
upper gradient of the function 
$2Cu$ is admissible for $\Modp(\Ga(\itoverline{B},F_n))$.
Taking infimum over all such upper gradients  implies that
\[
 \Modp(\Ga(\itoverline{B},F_n)) 
\le \int_X g_{2Cu}^p\, d\mu = (2C)^p \cpDp(F_n,G_n).  
\]
On the other hand, if $\vint_{B} u\,d\mu\le\tfrac12$, then
applying the weak Harnack inequality~\eqref{eq-weak-Harnack-B_1} to
the \p-harmonic function $1-u$, we see that $2C(1-M)\ge1$.
Thus, any 
upper gradient of the function $2C(1-u)$ 
is admissible for $\Modp(\Ga(\itoverline{B},G_n))$ and hence
\[
 \Modp(\Ga(\itoverline{B},G_n)) 
\le \int_X g_{2C(1-u)}^p\, d\mu = (2C)^p \cpDp(F_n,G_n). 
\]
Combining the above two inequalities, we have
\[
\lim_{n\to\infty} \cpDp(F_n,G_n)
\ge \frac{1}{(2C)^p} \lim_{n\to\infty} 
\min \{\Modp(\Ga(\itoverline{B},F_n)),\Modp(\Ga(\itoverline{B},G_n))\}
>0.\qedhere
\]
\end{proof}

\begin{lem}
Suppose that $\{F_n\}_{n=1}^\infty$ and $\{G_n\}_{n=1}^\infty$ are two disjoint
decreasing sequences
of closed nonempty sets in $X$
satisfying~\ref{l-i} of 
Definition~\ref{def-hyp-seq}. Assume that
$\Modp(\Ga(F_1,G_1)) <\infty$ and that
\[
\lim_{n\to\infty}\Modp(\Ga(F_n,G_n)) =: 2c_0 >0.
\]
Then
$\{F_n\}_{n=1}^\infty$ and $\{G_n\}_{n=1}^\infty$ are \p-hyperbolic sequences.
\end{lem}

\begin{proof}
Since, by  Lemma~\ref{lem-Mod=capDp},
\[
2c_0 \le \cpDp(F_1,G_1)<\infty,
\] 
there is a function $u\in\Dp(X)$ such
that $0\le u\le 1$ on $X$, $u=0$ on $F_1$ and $u=1$ on $G_1$.
Let $B$ be a sufficiently large ball
such that 
\[
\int_{X\setm B} g_u^p\,d\mu < c_0. 
\]
By changing $g_u$ on a set of zero measure, we can assume that it is a
Borel function.

Let $n$ be large enough so that
$F_n$ and $G_n$ are disjoint from $\clB$.
The curve family $\Ga_n=\Ga(F_n,G_n)$
can be written as the union $\Ga'\cup \Ga''$,
where $\Ga'$ contains the curves from $\Ga_n$ passing through $\clB$,
while $\Ga''$ consists of 
those curves from $\Ga_n$ which avoid $\clB$.
By the choice of $u$, 
the function $\rho:= g_u \chi_{X\setm\clB}$ is 
admissible for $\Modp(\Ga'')$ and hence
$\Modp(\Ga'')<c_0$.
Since every curve in $\Ga'$ has a subcurve in $\Ga(\clB,F_n)$, it
follows from \cite[Lemma~1.34]{BBbook} or \cite[p.~128]{HKST} that
\[
\Modp(\Ga(\clB,F_n)) \ge \Modp(\Ga')\ge \Modp(\Ga_n) - \Modp(\Ga'') >c_0,
\]
and similarly $\Modp(\Ga(\clB,G_n))>c_0$.
Thus, by  Lemmas~\ref{lem-Mod=capDp} and~\ref{lem-cap-K1-K2},
$\{F_n\}_{n=1}^\infty$ and $\{G_n\}_{n=1}^\infty$ are \p-hyperbolic sequences.
\end{proof}

\begin{proof}[Proof of Theorem~\ref{thm-2-hyp-seq}]
  Since $\cpDp(F_1,G_1)=\Modp(\Ga(F_1,G_1))<\infty$ 
  by Lemma~\ref{lem-Mod=capDp}, 
there is $f \in \Dp(X)$ such that $f \equiv 0$
on $F_1$ and $f \equiv 1$ on $G_1$.

As in the proof of Lemma~\ref{lem:hyp-seq-imp-Ga},
for each $n\ge1$, let
$u_n=H_{X\setminus (F_n\cup G_n)} f$ be
  the \p-harmonic extension
  of $f$ in $X \setm (F_n \cup G_n)$.
By Proposition~\ref{prop-cpDp-Hf},
\[
\int_X g_{u_n}^p\, d\mu=\cpDp(F_n,G_n)\le \cpDp(F_1,G_1).
\]

Lemma~\ref{lem-reflex-conv} 
provides us with a \p-harmonic function $u$ on $X$ and convex combinations 
\[
v_n=\sum_{j=n}^{N_n}\la_{j,n}u_j, \quad \text{where $0\le \la_{j,n}\le 1$ and
$\sum_{j=n}^{N_n}\la_{j,n}=1$,}
\]
such that $v_n \to u$ locally uniformly and
$\|g_{v_n}-g_u\|_{L^p(X)}\to 0$ 
as $n\to\infty$.
Hence
\begin{align*}
\biggl( \int_X g_u^p\, d\mu \biggr)^{1/p} 
&= \lim_{n\to\infty}  \biggl( \int_X g_{v_n}^p\, d\mu \biggr)^{1/p} 
\le \lim_{n\to\infty} \sum_{j=n}^{N_n}  \la_{j,n} 
        \biggl( \int_X g_{u_j}^p\, d\mu \biggr)^{1/p} \\
&\le \cpDp(F_1,G_1)^{1/p}  <\infty,
\end{align*}
showing that $u$ is a bounded \p-harmonic function in $X$ with finite energy.

Moreover, each $v_n$ is admissible for 
$\cpDp(F_{N_n},G_{N_n})$, and so by Lemma~\ref{lem:hyp-seq-imp-Ga},
\[
 \lim_{n\to\infty} \int_X g_{v_n}^p\, d\mu 
\ge \lim_{n\to\infty} \cpDp(F_{N_n},G_{N_n}) > 0,
\]
showing that $\int_X g_u^p\, d\mu>0$ and so $u$ is nonconstant.
\end{proof}

\section{Classification of metric measure spaces} \label{sec:Classification}

Recall from Definition~\ref{def-qmin} and the equality $\Np_0(X)=\Np(X)$
that a function $u\in \Np\loc(X)$ is quasiharmonic
in $X$ if for each $\phi\in \Np(X)$
we have
\begin{equation} \label{eq-qmin-char}
\int_{\phi \ne 0} g_u^p\, d\mu \le Q_u\int_{\phi \ne 0} g_{u+\phi}^p\, d\mu.
\end{equation}

\begin{thm}\label{thm:HD=HBD}
If $X$ supports a nonconstant quasiharmonic function with finite
energy, then $X \notin O^p_{HBD}$.
In particular,
\[
O^p_{HD}=O^p_{HBD}=O^p_{QD}=O^p_{QBD}.
\]
\end{thm}

A direct consequence of this result together with Theorem~\ref{thm-Liouville}
  is the following improvement of one of the main results in
  Bj\"orn--Bj\"orn--Shan\-mu\-ga\-lin\-gam~\cite[Theorem~1.1]{BBS-Liouville}
  under the additional assumption that $X$ is complete
  (which under our standing
    assumptions follows from the properness of $X$).
  In Section~\ref{sect-noncomp} we explain how to obtain it for
  noncomplete spaces.
  
\begin{cor} \label{cor-finite-energy-Liouville}
 Assume that $\mu$ is globally doubling and supports a global
  \p-Poincar\'e inequality.
  If $u \in \Dp(X)$ is a quasiharmonic function on $X$ with
  finite energy, then it is constant. 
In particular, $X \in O^p_{QD}$.
\end{cor}

\begin{proof}[Proof of Theorem~\ref{thm:HD=HBD}]
We prove the contrapositive statement.
So assume that $X\in O^p_{HBD}$
and that $u$ is a quasiharmonic function on $X$
with finite energy 
$\int_X g_u^p\, d\mu<\infty$.
Our aim is to show that $u$ is constant, and we do so by showing that
$g_u=0$ a.e.\ in~$X$.

Fix $x_0 \in X$ and let $B_j=B(x_0,j)$, $j=1,2,\ldots$\,.
As $X$ is unbounded, $\Cp(X\setminus B_j)>0$.
For each positive integer $k$ let $u_k=\min\{k,\max\{-k,u\}\}$. 
Let $v_{k,j}=H_{B_j} u_k$ be
  the \p-harmonic extension
  of $u_k$ in $B_j$.
Then $|v_{k,j}|\le k$ and 
\[
\int_X g_{v_{k,j}}^p\, d\mu
  =\int_{B_j} g_{v_{k,j}}^p\, d\mu +\int_{X\setminus B_j} g_{u_k}^p\, d\mu
\le \int_X g_{u_k}^p\,d\mu 
< \infty.
\]
Lemma~\ref{lem-reflex-conv} provides us with 
convex combinations $\vhat_{k,j}$ of the sequence
$\{v_{k,j}\}_{j=1}^\infty$
which converge
locally uniformly in $X$ to a bounded function $v_k\in\Dp(X)$
that is \p-harmonic in $X$,
and moreover 
\[
\|g_{\vhat_{k,j}}-g_{v_k}\|_{L^p(X)}\to 0, \quad \text{as } j\to\infty.  
\]
As we have assumed that $X\in O^p_{HBD}$ (at the beginning of the proof),
$v_k$ must be constant on $X$.
Thus $g_{v_k}=0$ and $g_{\vhat_{k,j}} \to 0$ in $L^p(X)$ as $j\to\infty$.

Since $\phi_{k,j}:= \vhat_{k,j}-u_k$ are
convex combinations of functions in $\Np(X)$,
we see that  $\phi_{k,j} \in\Np(X)$
and
\(
g_{u+\phi_{k,j}} \le g_{u-u_k} + g_{\vhat_{k,j}}.
\)
The quasiminimizing property~\eqref{eq-qmin-char} of $u$ then implies that
\begin{align*}
\int_X g_u^p\,d\mu &= \int_{\phi_{k,j}\ne0}g_u^p\, d\mu+ \int_{\phi_{k,j}=0}g_u^p\, d\mu\\
&\le Q_u \int_{\phi_{k,j}\ne0} g_{u+\phi_{k,j}}^p \,d\mu
+ \int_{\phi_{k,j}=0} g_{u+\phi_{k,j}}^p \,d\mu \\
&\le Q_u \int_X ( g_{u-u_k} + g_{\vhat_{k,j}} )^p \,d\mu \\
&\le 2^p Q_u \biggl( \int_{|u|>k} g_u^p \,d\mu
+ \int_X g_{\vhat_{k,j}}^p \,d\mu \biggr),
\end{align*}
where $Q_u$ is the quasiminimizing constant associated with $u$.
Letting $j\to\infty$ and then $k\to\infty$ shows that
$g_u=0$  a.e.~in $X$.
From the local Poincar\'e inequality, 
the connectivity of $X$
and the continuity of $u$
we conclude that $u$ must be constant on~$X$.
\end{proof}

\begin{remark} \label{rmk-unifPI}
It follows directly from Theorem~\ref{thm:HD=HBD} that 
the following two equivalent conditions can be added to Theorem~10.5 in
Bj\"orn--Bj\"orn--Shan\-mu\-ga\-lin\-gam~\cite{BBSunifPI}:
\begin{enumerate}
  \addtocounter{enumi}{2}
  \item There exists a nonconstant bounded \p-harmonic function on $(X,d,\mu)$
with finite \p-energy.
  \item There exists a nonconstant quasiharmonic function on $(X,d,\mu)$
with finite \p-energy.
\end{enumerate}
Similar modifications can also be made in the conclusions in
\cite[Example~10.8]{BBSunifPI}.
\end{remark}

We are now ready to state and prove the converse of
  Theorem~\ref{thm-2-hyp-seq}.

\begin{prop}   \label{prop-p-harm-imp-2-hyp-seq}
If $X$ supports a nonconstant bounded \p-harmonic function with finite
energy, then there are two disjoint \p-hyperbolic sequences
$\{F_n\}_{n=1}^\infty$ and $\{G_n\}_{n=1}^\infty$ such that 
$\Modp(\Ga(F_1,G_1)) <\infty$.
In particular, $X$ is \p-hyperbolic.
\end{prop}

To prove Proposition~\ref{prop-p-harm-imp-2-hyp-seq}, we shall need the
following definition, which extends the well-known notion of
\p-parabolic spaces to open subsets, see
Proposition~\ref{prop-par-hyp}
below.
For manifolds and $p=2$, this definition
appeared in
Grigor$'$yan~\cite[Definition~3]{Grig87}, \cite[Section~14.1]{Grig99}
and for metric spaces and $p>1$ in
Hansevi~\cite[Definition~4.1]{hansevi2}.

\begin{deff}   \label{def-p-par-set}
An unbounded open set $\Om\subset X$ is \emph{\p-parabolic} if for each 
compact set $K\subset\Om$ there exist functions $u_j\in\Np(\Om)$
such that $u_j\ge1$ on $K$ for all $j=1,2,\ldots$ and
\begin{equation}   \label{eq-def-p-par-set}
\int_\Om g_{u_j}^p\,d\mu \to 0 \quad \text{as } j\to\infty.
\end{equation}
\end{deff}

  \begin{prop} \label{prop-par-hyp}
$X$ is \p-parabolic in the sense of Definition~\ref{def-p-par-set}
if and only if it is \p-parabolic in the sense of Definition~\ref{def-hyp-end}. 
  \end{prop}

\begin{proof}
If $X$ is \p-hyperbolic 
in the sense of
    Definition~\ref{def-hyp-end}, then
fixing $x_0\in X$ and $K:=\itoverline{B(x_0,1)}$,
we know that
\[
\lim_{n\to\infty}  \Modp(\Gamma(K,X\setminus B_n)) = :c >0,
\]
where $B_n=B(x_0,n)$.
Now suppose that there is a sequence $u_j\in N^{1,p}(X)$ as in
Definition~\ref{def-p-par-set},
related to the compact set $K$, and for each
positive integer $n>2$ let $\eta_n$ be a $1$-Lipschitz function on $X$ 
such that $0\le \eta_n\le1$ on $X$,
$\eta_n=1$ on $\clB_{n-1}$, and $\eta_n=0$ outside $B_{n}$.
Then $v_{n,j}:=\eta_nu_j\in N^{1,p}(X)$
with $v_{n,j}=1$ on the compact set $K$ and $v_{n,j}=0$ outside
$B_{n}$. 
It then follows from the definition of \p-weak upper gradients that
for \p-almost every curve $\gamma\in \Gamma(K,X\setminus B_n)$ 
we have that  $1\le \int_\gamma g_{v_{n,j}}\,ds$.
Since 
\[
g_{v_{n,j}}\le u_j\chi_{B_{n}\setminus B_{n-1}}+g_{u_j}\chi_{B_{n}},
\] 
we see that 
\begin{equation*}
\Modp(\Gamma(K,X\setminus B_n))
\le 2^{p} \biggl(\int_{X\setminus B_{n-1}}|u_j|^p\, d\mu
           +\int_{B_{n}}g_{u_j}^p\, d\mu\biggr).
\end{equation*}
Letting $n\to\infty$ gives us that
\[
0<c \le 2^p\int_X g_{u_j}^p\, d\mu,
\]
which then forbids the sequence $u_j$ from satisfying~\eqref{eq-def-p-par-set},
that is, $X$ cannot be \p-parabolic in the sense of Definition~\ref{def-p-par-set}.

Conversely, if $X$ is not \p-parabolic in the sense of
Definition~\ref{def-p-par-set}, then there exist $c_0$ and a compact
set $K_0\subset X$ such that for every $u\in\Np(X)$ with $u\ge1$ on
$K_0$,
\[
\int_X g_u^p\, d\mu \ge c_0 >0.
\]
In particular, 
\[
\lim_{n\to\infty}\cpDp(K_0, X\setminus B_n) \ge c_0,
\]
which in combination with Lemmas~\ref{lem-Mod=capDp}
and~\ref{lem-cap-K1-K2} implies that
\[
\lim_{n\to\infty}\Modp(\Gamma(\clB_1, X\setminus B_n))>0,
\]
that is, $X$ is \p-hyperbolic in the sense of
Definition~\ref{def-hyp-end}.
\end{proof}

\begin{proof}[Proof of Proposition~\ref{prop-p-harm-imp-2-hyp-seq}]
Suppose that $u\in\Dp(X)$  is a bounded nonconstant \p-harmonic function $u$ on $X$ 
with finite energy.
Without loss of generality we may assume that
\[
\inf_X u=-1 \quad \text{and} \quad \sup_X u=2.
\]
Setting $\Om=\{x: u(x)<0\}$, choose a point $x_0\in\Om$. 
For  $n=1,2,\ldots$\,, let $F_n=\clOm\setm B(x_0,n)$.
We shall show that the sequence $\{F_n\}_{n=1}^\infty$ is \p-hyperbolic.

Assume not. Let $K\subset\Om$ be an arbitrary compact set.
Then by Lemmas~\ref{lem-Mod=capDp} and~\ref{lem-cap-K1-K2},
$\cpDp(K,F_n)=\Modp(\Ga(K,F_n))\to 0$ as $n\to\infty$.
In particular, for sufficiently large $n$,
there exist
$u_n \in \Dp(X)$ such that $u_n=1$ in $K$, $u_n=0$ in $F_n$,
$0 \le u_n \le 1$ in $X$ and
\[
\int_X g_{u_n}^p \,d\mu \to 0
\quad \text{as } n \to \infty.
\]
Since 
$u_n|_\Om$ has bounded support we see that $u_n \in \Np(\Om)$.
As $K$ was arbitrary, we conclude that $\Om$ is \p-parabolic in the sense of
Definition~\ref{def-p-par-set}.
Since $u=0$ on $\bdy\Om$, applying Corollary~7.7 in Hansevi~\cite{hansevi2}
(see Remark~\ref{rem-local-ass-Hansevi})
to
the constant function $f\equiv0$ then implies that $u\equiv0$ in $\Om$, which
is a contradiction.
Thus, the sequence $\{F_n\}_{n=1}^\infty$ is \p-hyperbolic,
  and hence $X$ is \p-hyperbolic.

Similarly, considering $\Om'=\{x: u(x)>1\}$ and $x'_0\in\Om'$,
we conclude that $G_n=\clOm'\setm B(x'_0,n)$ also forms a \p-hyperbolic
sequence.
Clearly, the two sequences are disjoint.
Moreover, any upper gradient $g$ for $u$ is admissible for
$\Modp(\Ga(F_1,G_1))$ and hence $\Modp(\Ga(F_1,G_1))<\infty$.
\end{proof}

\begin{proof}[Proof of Theorem~\ref{thm-intro-hyp-seq}]
One implication follows
directly from Theorem~\ref{thm-2-hyp-seq}, while
the other (and the \p-hyperbolicity) follows
from Proposition~\ref{prop-p-harm-imp-2-hyp-seq}.
The last part (about \p-hyperbolic ends) follows from
Corollary~\ref{cor:mainly-here}.
\end{proof}

\begin{proof}[Proof of Theorem~\ref{thm-intro-class}]
The inclusions
\begin{equation*}
  \begin{matrix}
    O^p_{HP} & \subset & O^p_{HB} & \subset & O^p_{HBD}  \\
   \vsubset & & \vsubset && \vsubset \\
  O^p_{QP} & \subset &  O^p_{QB} & \subset & O^p_{QBD}
      \end{matrix}
\end{equation*}
are trivial.
That $O^p_{HD}=O^p_{HBD}=O^p_{QD}=O^p_{QBD}$ follows
from Theorem~\ref{thm:HD=HBD},
while the inclusion
$O^p_{\para} \subset O^p_{HBD}$
follows from Theorem~\ref{thm-intro-hyp-seq}.
We have thus shown all inclusions.

As (unweighted) $\R^n \in O^p_{QP} \subset O^p_{HBD}$
for all $1<p<\infty$ and $n \ge 1$
(e.g.\ by Theorem~\ref{thm-Liouville}),
but
is \p-parabolic only for $p\ge n$,
we see  that
$O^p_{\para} \subsetneq O^p_{HBD}$
and that $\R^n \in O^p_{QP} \setm O^p_{\para}$ if $p<n$.

By Example~\ref{ex-R-hyp-par} below, there is a measure
$\mu$ on $\R$ (satisfying our standing assumptions)
such that
\[
(\R,\mu)\in O^p_{QB} \setm O^p_{HP}.
\]
It follows directly that
\[
  (\R,\mu)\in O^p_{HB} \setm O^p_{HP},
\quad \text{and} \quad
(\R,\mu)\in O^p_{QB} \setm O^p_{QP}.
\]

Finally, consider
  the Poincar\'e $n$-ball $B_\alp^n$ as in
  Sario--Nakai--Wang--Chung~\cite[Section~I.2.4]{SNWC},
  namely $B_\alp^n= B(0,1) \subset \R^n$, $n \ge 3$, equipped
  with the Poincar\'e-type metric 
  $ds_{\al}= (1-|x|^2)^\al\,dx$, $\al \le -1$,
  and the corresponding Lebesgue measure.
  This makes $B_\alp^n$ into an unbounded proper Riemannian manifold
  (and thus metric space) satisfying our
  standing assumptions. 
  By \cite[Lemma~I.2.8 and~I.2.9]{SNWC},
  $X \in O^2_{HD}\setm O^2_{HB}$.
\end{proof}

\section{\texorpdfstring{$\Dp(X)=N^{1,p}(X)+\R$}{Dp(X)=N1p(X)+R}} \label{sect-Dirichlet}

This section is devoted to Theorem~\ref{thm:main2},
  and we start with its proof. The
rest of the section discusses 
the converse of Theorem~\ref{thm:main2}.

\begin{proof}[Proof of Theorem~\ref{thm:main2}]
Since $X\notin O^p_{HD}$, there is a nonconstant
\p-harmonic function $u\in \Dp(X)$.
Suppose that there is some $c\in\R$ such that  $u+c\in N^{1,p}(X)$. 
Then $u+c$ is also nonconstant and \p-harmonic on $X$,
but this is in contradiction with 
Lemma~\ref{lem-Liouville-Np}.
\end{proof}

The following example shows that
$\Dp(X)=N^{1,p}(X)+\R$ can fail even when $X\in O^p_{HD}=O^p_{QD}$.

\begin{example} \label{ex:OpHD-but-not-hom}
Let $X=\R^n$ (unweighted) with $p > n \ge 1$
and let
\[
    u(x)= \sum_{j=0}^{\infty} \bigl(1-2^{-j}|x-(4^j,0,\ldots,0)|\bigr)_\limplus.
\]
Then both $\{x:u(x)=0\}$ and $\{x:u(x)>\tfrac12\}$ have infinite
measure and thus $u \notin N^{1,p}(X)+\R$.
However,
\[
  \int_{\R^n} g_u^p \, dx
  = \sum_{j=0}^\infty  2^{j(n-p)} \omega_n  < \infty,
\]
where $\om_n$ is the volume of the unit ball in $\R^n$.
Thus $u \in \Dp(X)$.

Note that $X \in O^p_{QD}$ by Corollary~\ref{cor-finite-energy-Liouville},
and even  $X \in O^p_{QP} $ 
by Theorem~\ref{thm-Liouville}.
\end{example}  

In Proposition~\ref{prop-OpHD-andHom} below
we show that the converse of Theorem~\ref{thm:main2} holds
provided that $X$ supports the following
global $(p,p)$-Sobolev inequality.

\begin{deff}\label{def:pp-Sobolev}
$X$ supports a \emph{global $(p,p)$-Sobolev inequality}
if there is a constant $C>0$ such that 
\begin{equation} \label{eq-pp-Sob}
\int_X|u|^{p}\, d\mu\le C \int_Xg_u^p\, d\mu
\end{equation}
whenever $u\in N^{1,p}(X)$.
\end{deff}

One can equivalently
require  \eqref{eq-pp-Sob} to just hold
for bounded $u \in\Np(X)$ with bounded support,
see \cite[Proposition~7.1.35]{HKST}.
  If $X$ is a simply connected complete Riemannian
  manifold with sectional curvature $K  \le -a^2 <0$,
  then it supports a global $(p,p)$-Poincar\'e inequality,
  see
Holopainen--Lang--V\"ah\"akangas~\cite[p.~129]{HLV}.

The global $(p,p)$-Sobolev inequality holds
if and only if 
the \emph{Rayleigh quotient} 
\[
R_p(X):=\inf_u  \frac{\int_X g_u^p\, d\mu}{\int_X|u|^{p}\, d\mu}>0,
\]
where the infimum is taken over all $u \in \Np(X)$ with $\|u\|_{\Np(X)} > 0$.

Classically, the Rayleigh quotient for $p=2$ equals the first
eigenvalue $\la_1$
(the bottom of the spectrum) of the
$2$-Laplacian.
In the nonlinear case, the Rayleigh quotient
is associated with a nonlinear eigenvalue problem
that has even been studied on metric spaces (with upper gradients
as here),
see e.g.\ Garc\'\i a Azorero--Peral Alonso~\cite{GP} (on $\R^n$)
and
Latvala--Marola--Pere~\cite{LMP}.

It was shown in
Li--Wang~\cite[Theorem~1.4\,(2)]{LW} that if $M$ is a complete $2$-hyperbolic Riemannian manifold 
with $\lambda_1>0$, then the measure of balls centered at a point in $M$
grows at least exponentially with 
respect to the radius.
A similar exponential volume growth was identified in
Buckley--Koskela~\cite[Theorem~0.1\,(2)]{BuK} for 
proper \p-hyperbolic metric measure spaces supporting a
global $(p,p)$-Sobolev inequality.

\begin{prop} \label{prop-OpHD-andHom}
If in addition to our standing assumptions on $X$, we know that
$X\in O^p_{HD}$ and $X$ supports a global $(p,p)$-Sobolev inequality,
then 
\[
\Dp(X)=N^{1,p}(X)+\R.
\]
\end{prop}  

\begin{proof}
Let $f \in \Dp(X)$, $x_0 \in X$ and set $B_j=B(x_0,j)$, $j=1,2,\ldots$\,.
As $X$ is unbounded, $\Cp(X\setminus B_j)>0$.
For each positive integer $k$, let 
$v_{k}=H_{B_k} f$ be
  the \p-harmonic extension
  of $f$ in $B_k$.
By the global $(p,p)$-Sobolev inequality, 
\begin{align}
\int_X|f-v_k|^{p}\, d\mu
&\le C\int_Xg_{f-v_k}^p\, d\mu \nonumber \\
& \le 2^{p}C\biggl(\int_Xg_f^p\, d\mu+ \int_Xg_{v_k}^p\, d\mu\biggr)
\le 2^{p+1}C\int_Xg_f^p\, d\mu. 
\label{eq-f-vk}
\end{align}
Thus for each $j=1,2,\ldots$\,, 
\begin{align*}
   \int_{B_j}|v_k|^p\, d\mu
   & \le 2^p \biggl(\int_{B_j}|f-v_k|^p\, d\mu + \int_{B_j}|f|^p\, d\mu\biggr) \\
   & \le 2^{2p+1} C\biggl(\int_{X} g_f^p\, d\mu
        + \int_{B_j}|f|^p\, d\mu\biggr). \\
\end{align*}
As $\Nploc(X)=\Dploc(X)$ 
(see Section~\ref{sect-prelim}),
and $X$ is proper, we see that $f \in \Np(B_j) \subset L^p(B_j)$,
and so $\{v_k\}_{k=1}^\infty$ is a bounded sequence in $N^{1,p}(B_j)$. 

Now an argument as in the proof of Lemma~\ref{lem-reflex-conv}
shows that we 
have a sequence $\{\hat{v}_k\}_{k=1}^\infty$, of convex combinations of
$\{v_k\}_{k=1}^\infty$, that converges in $N^{1,p}(B_j)$ 
for each $j$
to a 
function $v$ on $X$, with $g_v\in L^p(X)$.
By Shan\-mu\-ga\-lin\-gam~\cite[Theorem~1.1]{Sh03}, the function $v$ is
  \p-harmonic in each $B_j$ and thus in $X$.
Since $X\in O^p_{HD}$, the function $v$ must be constant, say $v \equiv c$.
As in \eqref{eq-f-vk},
\[
\int_{B_j}|f-c|^p\, d\mu 
= \lim_{k \to \infty} \int_{B_j}|f-\hat{v}_k|^p\, d\mu 
\le \limsup_{k \to \infty} \int_X|f-\hat{v}_k|^p\, d\mu
\le C' \int_Xg_f^p\, d\mu.
\]
After letting also $j\to \infty$, 
we see that $\int_X|f-c|^p\, d\mu$ is finite, that is, $f-c\in N^{1,p}(X)$.
The converse inclusion is trivial.
\end{proof}

\section{Examples related to Liouville type classes} \label{sec-ex}

In this section we explore some examples that illustrate the (non)equality of the
Liouville classes.
As mentioned in the introduction, the Euclidean space $\R^n$ has
only one end when $n\ge 2$, and that end
is \p-hyperbolic only when $1\le p<n$. On the other hand, $\R$ has two distinct ends.

We begin with an example showing that $\R^2$
  can be equipped with a weight so that it has two well-separated
  \p-hyperbolic sequences, even though
  it only has one end, cf.\ Theorem~\ref{thm-intro-hyp-seq}.

\begin{example}  \label{ex-weighted-R2}
Let $X=\R^2$ be equipped with the Euclidean distance and the measure
$d\mu=w\,dx$, 
where
\[
w(x)=e^{-\dist(x,A)}
\quad \text{and} \quad  
A=\{x=(x_1,x_2):|x_2| \le |x_1|\}.
\]
Also let $1<p<2$.

Observe that as $w$ is ``uniformly almost constant'' on every ball of radius $1$,
$\mu$ is uniformly locally doubling and supports a uniformly local $1$-Poincar\'e
inequality.

Even though $X$ only has one end, it
is still possible to use Theorem~\ref{thm-2-hyp-seq}
to show that $X$ supports a 
nonconstant bounded \p-harmonic function with finite energy.
Let $x_0=(0,0)$, $B=B(x_0,1)$ and
\begin{align*}
  F_n &=\{x=(x_1,x_2) \in A : x_1 \ge 2n \}, \\
  G_n &=\{x=(x_1,x_2) \in A : x_1 \le -2n \}, 
\quad n=1,2,\ldots.
\end{align*}
Now, by symmetry, 
\begin{align}   \label{est-hyp-Fn-hyp-R2}
 \Modp(\Ga(\itoverline{B},\R^2 \setm (-2n,2n)^2))
 & \ge \Modp(\Ga(\itoverline{B},F_n)) \nonumber \\
& \ge\tfrac{1}{4} \Modp^{\R^2}(\Ga(\itoverline{B},\R^2 \setm (-2n,2n)^2)), 
\end{align}
where 
$\Modp^{\R^2}$ denotes the standard \p-modulus in unweighted $\R^2$ with
respect to the Lebesgue measure.
Since $1<p<2$, we know that unweighted $\R^2$ is \p-hyperbolic and hence
the right-hand side in \eqref{est-hyp-Fn-hyp-R2} has a
  positive lower bound as $n\to\infty$.
It follows that $X$ is \p-hyperbolic and that
$\{F_n\}_{n=1}^\infty$ is a \p-hyperbolic sequence in $X$.
Similarly $\{G_n\}_{n=1}^\infty$ is a \p-hyperbolic sequence in $X$.

Next, every curve connecting $F_1$ to $G_1$ must pass through 
the strip 
\[
S:= \{x\in\R^2: |x_1|\le\tfrac12\},
\] 
whose characteristic function $\chi_S$ is thus admissible 
for the family $\Ga(F_1,G_1)$ 
of such curves.
A simple calculation then shows that
\[
\Modp(\Ga(F_1,G_1)) \le \int_{\R^2} \chi_S \, d\mu < \infty.
\]
Hence  Theorem~\ref{thm-2-hyp-seq} is applicable and provides us with
a bounded \p-harmonic function in $X$ with finite energy.

Note that $\R^2$ equipped with the Lebesgue measure does not
 support any nonconstant bounded \p-harmonic function.
  It therefore follows  from Proposition~\ref{prop-p-harm-imp-2-hyp-seq}
  that  $\Modp^{\R^2}(\Ga(F_n,G_n))=\infty$ for each positive integer~$n$.
\end{example}

We have seen that spaces with at least two \p-hyperbolic ends
support nonconstant bounded \p-harmonic functions with finite energy,
while parabolic spaces do not.
For spaces with only one end, which is \p-hyperbolic,
Example~\ref{ex-weighted-R2} and unweighted
$\R^n$ with $1<p<n$ show that they may
or may not support nonconstant bounded \p-harmonic functions with finite energy.
A natural question is what happens in spaces with only 
one
\p-hyperbolic end and at least one \p-parabolic end.
The following examples and Proposition~\ref{prop-hyp+par-end} show
that both situations are possible.

\begin{example}  \label{ex-hyp+par-end-R2}
Consider  
\[
X= \{x=(x_1,x_2)\in\R^2: x_2\le0\} \cup \bigl([-1,1] \times (0,\infty)\bigr),
\]
equipped with the Euclidean distance and the measure
$e^{-\dist(x,A)}\,dx$, where
\[
A=\{(x_1,x_2): -|x_1|\le  x_2\le0\} \cup \bigl([-1,1] \times (0,\infty)\bigr).
\]
Also let $1<p<2$.
Similar to Example~\ref{ex-weighted-R2}, we see that $X$ has one
\p-hyperbolic end  and contains two disjoint \p-hyperbolic sequences
\[
F_n =\{x=(x_1,x_2) \in A : x_1\ge 2n \} \quad \text{and} \quad
G_n =\{x=(x_1,x_2) \in A : x_1\le -2n \},
\]
$n=1,2,\ldots$\,,
while the strip $[-1,1] \times (0,\infty)$ forms a \p-parabolic end.
The uniformly local doubling property and a uniformly local
1-Poincar\'e inequality are also satisfied.
Theorem~\ref{thm-2-hyp-seq} now implies that $X$ supports a
nonconstant bounded \p-harmonic function with finite energy.
\end{example}

In the following example
we will see that when suitably equipped with a weight, $\R$
carries a nonconstant positive $2$-harmonic function
but no nonconstant bounded $2$-harmonic function.

\begin{example} \label{ex-R-hyp-par}
Consider $\R$ 
equipped with the Euclidean distance and the weight 
\[
w(x)= \begin{cases} 
           |x|^\alpha, &x\le-1, \\
           1, &x\ge-1,  \end{cases}
\]
for some fixed $\al>-1$. 
The measure
$d\mu(x)=w(x)\, dx$ is uniformly locally
doubling and supports a uniform local $1$-Poincar\'e inequality.
As pointed out above, $\R$ has two ends, denoted $\infty$ and $-\infty$.
By Proposition~\ref{prop-hyp+par-end} below, the end at $\infty$ is
\p-parabolic for each $p>1$, while the end at $-\infty$ is
\p-hyperbolic if (and only if) $1<p<1+\alpha$
(which then also requires $\alpha>0$). 
Moreover, $(\R,\mu)$ is a \p-hyperbolic space
in this case, since it   has a \p-hyperbolic end.
It thus follows, from Proposition~\ref{prop-hyp+par-end} again,
that
\[
(\R,\mu) \in (O^p_{QB}\cap O^p_{QD})\setminus
(O^p_{HP} \cup O^p_{\para})
\quad \text{when }1<p<1+\alpha.
\]
\end{example}

\begin{prop}  \label{prop-hyp+par-end}
Consider the real line $\R$, equipped with the Euclidean distance and 
the measure $d\mu=w\,dx$,
where $\mu$ is  locally doubling and supports a local \p-Poincar\'e inequality.
Then the following are true.
\begin{enumerate}
\item \label{f-1}
Each quasiharmonic function {\rm(}with respect to $\mu${\rm)}
on $\R$ is bounded if and only if it has finite energy.
\item \label{f-2}
The end at $\infty$ is \p-hyperbolic if and only if 
\begin{equation}  \label{eq-assume-w-hyp-pos}
\int_0^\infty w^{1/(1-p)}\,dx <\infty.
\end{equation}
\item \label{f-3}
The end at $-\infty$ is \p-hyperbolic if and only if 
\begin{equation}  \label{eq-assume-w-hyp-neg}
\int_{-\infty}^0 w^{1/(1-p)}\,dx <\infty.
\end{equation}
\item \label{f-4}
  The space $(\R,\mu) \in O^p_{\para}$ 
  if and only if both \eqref{eq-assume-w-hyp-pos} and \eqref{eq-assume-w-hyp-neg}
  fail.
\item If both \eqref{eq-assume-w-hyp-pos} and \eqref{eq-assume-w-hyp-neg} hold
  then there exists a nonconstant bounded global \p-harmonic function
  with finite energy, i.e.\ $(\R,\mu) \notin O^p_{HBD}$.
\item 
  If \eqref{eq-assume-w-hyp-pos} holds and \eqref{eq-assume-w-hyp-neg} fails\/
\textup{(}or \eqref{eq-assume-w-hyp-pos} fails and \eqref{eq-assume-w-hyp-neg}
  holds\/\textup{)},
  then
\[
(\R,\mu)\in (O^p_{QB} \cap O^p_{QD}) \setm O^p_{HP}.
\]
\item If both \eqref{eq-assume-w-hyp-pos} and \eqref{eq-assume-w-hyp-neg} fail
then $(\R,\mu)\in O^p_{QP} \cap O^p_{QD}$.
\end{enumerate}
\end{prop}

Weights $\mu$ on $\R$ as above were characterized in
Bj\"orn--Bj\"orn--Shan\-mu\-ga\-lin\-gam~\cite[Theorem~1.2]{BBSpadm}.
In particular it was shown that 
for each bounded interval $I$ there is a (global) $A_p$ weight $\wt$ on $\R$
such that $\wt=w$ on $I$.

\begin{proof}
\ref{f-1}
This follows from
\cite[Proposition~6.5]{BBS-Liouville}.

\ref{f-2}
To see that the end at $\infty$ is \p-hyperbolic when
\eqref{eq-assume-w-hyp-pos} holds, consider the
family $\Ga_R=\Ga([-1,0],[R,\infty))$, $R>0$,
and let $\rho$ be admissible for $\Modp(\Ga_R)$.
Since
\[
1 \le \int_0^{R} \rho\,dx 
   \le \biggl( \int_0^{R} \rho^p\,d\mu \biggr)^{1/p} 
            \biggl( \int_0^{R} w^{1/(1-p)}\,dx \biggr)^{1/p}  \le C \|\rho\|_{L^p(\R,\mu)},
\]
with $C$ independent of $R$,
letting $R\to\infty$ shows that the end at $\infty$ is \p-hyperbolic.
On the other hand, if \eqref{eq-assume-w-hyp-pos} fails then the function
\[
\rho_R(t) := \frac{w^{1/(1-p)}(t) \chi_{[0,R]}(t)}{\int_0^R w^{1/(1-p)}\,dx}
\]
is admissible for $\Modp(\Ga_R)$
(as we may assume that $w$ is a Borel function), with
\[
\int_\R \rho^p d\mu = \biggl( \int_0^R w^{1/(1-p)}\,dx \biggr)^{1-p}
  \to 0, \quad \text{as }  R\to\infty,
\]
and so the end at $\infty$ is \p-parabolic.
Thus \ref{f-2} has been shown, 
\ref{f-3} is shown similarly, and \ref{f-4} follows immediately.
 
The remaining statements follow from
\cite[Theorems~1.2 and~1.3]{BBS-Liouville}.
\end{proof}

As the next example shows, when a metric space has infinitely many ends, the \p-hyperbolicity of the space does
not imply the existence of a \p-hyperbolic end.

\begin{example}\label{ex-Tree-hyp-but-parEnds}
  The weighted infinite 
  binary tree $X$
from~\cite[Example~7.2]{BBS-Liouville} is an example of a space that
does not belong to $O^p_{HBD}$ for any $1<p<\infty$.
By Theorem~\ref{thm-intro-hyp-seq}, $X$ is \p-hyperbolic.
It is equipped with the geodesic metric, giving each edge unit length.
  Each geodesic ray, emanating from the root $v_0$, defines an end at infinity
  and corresponds to exactly one point in the so-called visual boundary of $X$.

  Fixing one such geodesic ray $\ga$ from the root, the measure on $X$
  is comparable to $2^{-\pi_\ga(x)}\,dm(x)$, where $m$ is the
  1-dimensional Lebesgue measure on each edge and
  $\pi_\ga(x)$ is the closest point on $\ga$ to $x$.
  Since the weight $w(x)=2^{-\pi_\ga(x)}$ is nonincreasing along each geodesic ray,
  an argument as in
Proposition~\ref{prop-hyp+par-end}\,\ref{f-2}
  together with Remark~\ref{rem:end-paths} tells us that the corresponding end
  must be \p-parabolic.
Thus, $X$ is a \p-hyperbolic space having only \p-parabolic ends.
It is also uniformly locally doubling
and supports a uniformly local $1$-Poincar\'e inequality.
\end{example}

\section{The finite-energy Liouville theorem in noncomplete spaces}
\label{sect-noncomp}

Recall that the standing assumptions from Section~\ref{sect-qmin}
are not required
in this section.

\begin{thm} \label{thm-Dp-Liouville-noncomplete}
  Assume that $X$ is a\/ \textup{(}not necessarily complete\/\textup{)}
  metric space
equipped with a globally doubling measure $\mu$ supporting
  a global \p-Poincar\'e inequality, where $1<p<\infty$.

  If $u \in \Dp(X)$ is a quasiharmonic function on $X$ with
  finite energy, then it is constant. 
\end{thm}

This shows that Theorem~1.1 in 
Bj\"orn--Bj\"orn--Shan\-mu\-ga\-lin\-gam~\cite{BBS-Liouville}
holds even if none of the sufficient conditions (a)--(d) therein is satisfied.

As in Definition~\ref{def-qmin}, a
function $u \in \Nploc(X)$ is a \emph{quasiminimizer} on the entire
  space $X$ if
\begin{equation*} 
     \int_{\phi \ne 0} g^p_u \,d\mu \le Q_u  \int_{\phi \ne 0} g^p_{u+\phi} \,d\mu
\quad \text{for all }\phi \in \Np(X),
\end{equation*} 
and a \emph{quasiharmonic function} is a continuous quasiminimizer.
However, the definition of quasiminimizers on \emph{strict subsets}
of noncomplete
spaces is more involved, 
see \cite[Section~3]{BBS-Liouville} for such a definition
and further discussion.

\begin{proof}
If $X$ is bounded, then this follows directly from
\cite[Theorem~1.1]{BBS-Liouville}.
So  we may assume that $X$ is unbounded.
As in \cite{BBnoncomp} we let $\Xhat$ be the completion
of $X$.
The metric $d$ extends directly to $\Xhat$
and we   define the complete Borel regular measure $\muhat$ on $\Xhat$ by letting
\[
     \muhat(E)=\mu(E \cap X)
\quad \text{for every Borel set } E \subset \Xhat,
\]
and then complete it, see \cite[Corrigendum]{BBnoncomp}.
It follows from \cite[Propositions~3.3 and~3.6]{BBnoncomp}
that $\muhat$ is globally doubling and supports
a global \p-Poincar\'e inequality on~$\Xhat$,
with the same doubling and Poincar\'e constants.

By \cite[Theorem~4.1]{BBnoncomp}, there is $\uhat \in \Dp(\Xhat)$
such that $\uhat=u$ $\CpX$-q.e.\ in $X$
and
\begin{equation} \label{eq-A0}
  g_{\uhat,\Xhat} \le A_0 g_{u,X}
  \quad \text{a.e.~in }X,
\end{equation}
where $A_0$ only depends on $p$, the global doubling 
constant and both constants in the global \p-Poincar\'e inequality.
Let $\phihat \in \Np(\Xhat)$ and $\phi=\phihat|_X$.
  Then,  $g_{u+\phi,X} \le g_{\uhat+\phihat,\Xhat}$
  a.e.\ in $X$, and thus using~\eqref{eq-A0},\
\[
   \int_{\phihat \ne 0}   g_{\uhat,\Xhat}^p \, d\muhat
   \le A_0^p    \int_{\phi \ne 0}   g_{u,X}^p \, d\mu
   \le A_0^p Q_u    \int_{\phi \ne 0}   g_{u+\phi,X}^p \, d\mu
   \le A_0^p Q_u    \int_{\phihat \ne 0}   g_{\uhat+\phihat,X}^p \, d\mu.
\]
Therefore
$\uhat$ is a quasiminimizer on $\Xhat$,
and hence 
has a continuous representative (see Section~\ref{sect-qmin})
that we can also call $\uhat$.
By \eqref{eq-A0}, we see that $\uhat$ has finite energy in $\Xhat$,
and thus by Corollary~\ref{cor-finite-energy-Liouville}, $\uhat$
is constant. As $u$ is continuous, it must also be constant.
\end{proof}

\end{document}